\documentclass[11pt,reqno]{elsarticle}

\makeatletter
\def\ps@pprintTitle{%
 \let\@oddhead\@empty
 \let\@evenhead\@empty
 \def\@oddfoot{}%
 \let\@evenfoot\@oddfoot}
\makeatother

\usepackage{amsmath,amsthm,amscd,wrapfig,amsfonts,amssymb,mathtools,siunitx}
\usepackage{enumitem,comment,mathabx,cool}
\usepackage{booktabs,multirow,lscape}
\usepackage{url,parskip,caption,breqn}

\usepackage[top=22mm, bottom=22mm, left=22mm, right=22mm]{geometry}

\newtheorem{theorem}{Theorem}

\newtheorem{remark}{Remark}

\DeclareMathOperator{\erf}{erf}

\providecommand{\Sr}{S(\mathbb{R})}

\providecommand{\e}[1]{\ensuremath{\times 10^{#1}}}

\makeatletter
\g@addto@macro\normalsize{%
  \setlength\abovedisplayskip{.4em}
  \setlength\belowdisplayskip{.4em}
  \setlength\abovedisplayshortskip{.4em}
  \setlength\belowdisplayshortskip{.4em}
}

\begin{document}

\begin{frontmatter}
\begin{abstract}
We describe a method for calculating the roots of 
special functions satisfying second order linear ordinary differential equations.
It exploits the recent observation that
the solutions of a large class of such equations 
can be represented  via nonoscillatory
phase functions, even in the high-frequency regime.
Our algorithm achieves near machine precision accuracy and
the time required to compute 
one   root of a solution
is independent of the frequency of oscillations of that
solution.     
Moreover, despite its great generality,
our approach is competitive with 
specialized, state-of-the-art methods for the construction of
Gaussian quadrature rules of large orders when it used in such
a capacity.
The performance of the scheme  is illustrated with several numerical experiments
and a Fortran implementation of our algorithm is available at the author's website.


\end{abstract}


\title
{
On the numerical calculation of the roots of special functions 
satisfying second order ordinary differential equations
}

\author[jb]{James Bremer}
\cortext[cor1]{Corresponding author (bremer@math.ucdavis.edu)}
\address[jb]{Department of Mathematics, University of California, Davis}

\end{frontmatter}

Special functions satisfying second order differential equations 
\begin{equation}
y''(t) + q(t) y(t) = 0
\ \ \ \mbox{for all}\ \ 
-\infty < a \leq t \leq b < \infty,
\label{introduction:second_order}
\end{equation}
where $q$ is smooth and positive,
are ubiquitous in mathematical physics,  and their roots 
play a number of roles.  Among other things,  
they are related to the resonances in mechanical and electromagnetic systems, 
arise in quantum mechanical calculations, and are the  nodes of Gaussian quadrature formulas.

Their role as the nodes of Gaussian quadrature formulas motivates
much of the interest in the numerical computation of the roots of special
function defined by equations of  the form (\ref{introduction:second_order}).
Every family of classical orthogonal polynomials --- 
Legendre polynomials, Hermite polynomials, Laguerre polynomials, etc. ---
is associated with a nonnegative weight function $w$ and a
collection of Gaussian quadrature rules, one for each positive integer $n$,
of the form
\begin{equation}
\int_a^b \varphi(t) w(t)\ dt  \approx
\sum_{j=1}^n \varphi(t_j) w_j.
\label{introduction:gaussian}
\end{equation}
The Gaussian quadrature rule (\ref{introduction:gaussian}) 
is  exact when
$\varphi$ is a polynomial of degree less than $2n$. 
The nodes $t_1,\ldots,t_n$  are the roots of
a polynomial $p$ of degree $n$ that satisfies
an equation of the form (\ref{introduction:second_order}),
and, at least in the case of the classical orthogonal polynomials,
 the weights  $w_1,\ldots,w_n$ can be calculated from the values 
of the derivatives of $p$ at the nodes $t_1,\ldots,t_n$.

Many schemes for the numerical computation of the nodes 
$t_1,\ldots,t_n$ and weights $w_1,\ldots,w_n$ of the rule (\ref{introduction:gaussian})
are based on the observation that classical orthogonal polynomials
also satisfy three-term recurrence relations.
Using such a recurrence relation to compute Newton iterates which converge to 
 the roots of an orthogonal polynomial
yields an $\mathcal{O}(n^2)$ method for the computation
of an $n$-point Gaussian quadrature rule.
The Golub-Welsch algorithm \cite{Golub-Welsch}
exploits the connection between three-term recurrence relations and the eigenvalues
of symmetric tridiagonal matrices in order to calculate the nodes
of an $n$-point Gaussian quadrature formula in $\mathcal{O}(n \log(n))$ operations;
it requires $\mathcal{O}\left(n^2\right)$ operations 
in order to compute both the nodes and weights.

In the last decade, several
$\mathcal{O}(n)$ methods for the calculation of $n$-point Gaussian quadrature 
rules have been proposed.  In \cite{Bogaert-Michiels-Fostier},  Newton's method
is combined with  a scheme, based on asymptotic formulas, 
for evaluating Legendre polynomials of arbitrary orders and arguments
in $\mathcal{O}(1)$ operations. 
A similar approach is taken in \cite{Hale-Townsend}, in which
asymptotic formulas for Jacobi polynomials are used to evaluate Newton iterates
which converge to the nodes of Gauss-Jacobi quadrature rules.
In \cite{Townsend-Trogdon-Olver}, the approach of  \cite{Hale-Townsend} is
applied to construct
Gauss-Hermite quadrature rules and 
 Gaussian quadrature rules for weight functions of the form
 $\exp(-V(x))$, where $V$ is a polynomial which grows at infinity.  
The  asymptotic formulas used in the construction of the later formulas
 are derived using  Riemann-Hilbert techniques.
Asymptotic formulas
which approximate the nodes and weights of 
Gauss-Legendre quadrature rules with double precision accuracy
are developed in \cite{Bogaert}.
Newton's method can be used to refine these approximations if greater accuracy
is required, but, in most cases, the asymptotic formulas
of \cite{Bogaert}  allow one to dispense with iterative methods entirely.
  The schemes of 
\cite{Bogaert,Bogaert-Michiels-Fostier,Hale-Townsend,Townsend-Trogdon-Olver,Segura2}
 all  have the property that any particular quadrature
node and its corresponding
weight can be calculated independently of the others,
making them suitable for parallelization.

The Glaser-Liu-Rokhlin method  \cite{Glaser-Rokhlin}
combines the Pr\"ufer transform with the classical Taylor series
method for the solution of ordinary differential equations.
It computes $n$ roots in  $\mathcal{O}(n)$ operations 
and is more general than the 
schemes  \cite{Bogaert, Bogaert-Michiels-Fostier, Hale-Townsend}
in that  it applies to special functions defined by second order differential equations
of the form
\begin{equation}
r_0(t) y''(t) + p_0(t) y'(t) + q_0(t) y(t) = 0,
\label{introduction:generalode}
\end{equation}
where $p_0$, $q_0$ and $r_0$ are polynomials of degree less than or equal to $2$.
This class includes the classical orthogonal polynomials, Bessel functions, 
prolate spheroidal wave functions, etc.
The Glaser-Liu-Rokhlin algorithm does, however, suffer from several disadvantages.
It is typically slower than the methods discussed above,
and is unsuitable for parallel implementation since
the roots must be computed in sequential order
(that is, the 
the computation of the $(n+1)^{st}$ root can only proceed once the $n^{th}$ root has been
obtained).    Moreover,  some precision is lost when the Glaser-Liu-Rokhlin
 technique is used to compute the 
weights of a Gaussian quadrature rule of large order
(see \cite{Bogaert-Michiels-Fostier} and \cite{Hale-Townsend} for discussions of this issue).
In \cite{Segura2},  a method for the calculation of  roots 
of solutions of quite general second order ordinary differential equations of the
form (\ref{introduction:second_order}).
It operates via a fixed point method which, unlike Newton's method, is guaranteed to converge.
It has the disadvantage, though, that the regions of monotonicity of the coefficient
$q$ must be explicitly known.

Here, we describe a fast and highly accurate algorithm for calculating the roots of special
functions which 
applies in even greater generality than the Glaser-Liu-Rokhlin method.
Indeed, it can be used to compute the roots of a solution of
an equation of the form (\ref{introduction:second_order})
as long as $q$ is a nonoscillatory function which
is positive and analytic on the interior of the interval
$[a,b]$ ($q$ can have poles, branch cuts or zeros at the endpoints $a$ and $b$).
Of course,   the ostensibly more general second order linear ordinary differential equation
(\ref{introduction:generalode})
can be easily  transformed into the form (\ref{introduction:second_order})
(see, for instance, \cite{Coddington-Levinson}),
and assuming that the coefficients $r_0$, $p_0$ and $q_0$ are nonoscillatory,
the algorithm of this paper can be brought to bear on the resulting equation.
Despite its great generality, 
our approach is competitive with the specialized
algorithm of \cite{Bogaert} 
for the computation of Gauss-Legendre quadrature rules of large orders,
and considerably faster than the specialized approach of \cite{Hale-Townsend}
for the computation of Gauss-Jacobi quadrature rules of large orders.
See the experiments of Sections~\ref{section:experiments:legendre}
and (\ref{section:experiments:jacobi}) for timings.  See also Section~\ref{section:experiments:laguerre},
where the algorithm of this paper is used to compute generalized Gauss-Laguerre quadrature rules.


Our approach exploits the fact
that when the coefficient $q$ in (\ref{introduction:second_order})
is nonoscillatory,
solutions of the equation (\ref{introduction:second_order})
can be represented to high accuracy via a nonoscillatory phase function, even
when  the magnitude of $q$ is large.
A smooth function
$\alpha:[a,b] \to \mathbb{R}$  is a phase function for
Equation~(\ref{introduction:second_order})
if $\alpha'$ is positive on $[a,b]$ and 
the pair of functions $u,v$ defined by the formulas
\begin{equation}
u(t) = \frac{\cos(\alpha(t)) }{\sqrt{\alpha'(t)}}
\label{introduction:u}
\end{equation}
and
\begin{equation}
v(t) = \frac{\sin(\alpha(t)) }{\sqrt{\alpha'(t)}}
\label{introduction:v}
\end{equation}
form a basis in the space of solutions of (\ref{introduction:second_order}).
Phase functions play a major role in the theories of special functions and
global transformations of ordinary differential equations
\cite{Boruvka,Neuman,NISTHandbook,Andrews-Askey-Roy}, and
are the basis of many numerical algorithms for the evaluation of special functions
(see \cite{Spigler-Vianello, Goldstein-Thaler,Heitman-Bremer-Rokhlin-Vioreanu}
for  representative examples).  

When $q$ is large in magnitude, most phase functions for (\ref{introduction:second_order})
are highly oscillatory.
  However, it has long been known that certain second order differential  
equations --- such Bessel's equation and Chebyshev's equation ---
  admit nonoscillatory phase functions.
In \cite{Heitman-Bremer-Rokhlin} and \cite{Bremer-Rokhlin}, it is shown
that, in fact, essentially all equations of the form (\ref{introduction:second_order}),
where $q$ is positive and nonoscillatory, admit a nonoscillatory 
 phase function $\alpha$ 
which represents solutions of (\ref{introduction:second_order}) with  high
accuracy.  
The function $\alpha$ is nonoscillatory in the sense that  it can be represented
using various series expansions (e.g., expansions in Chebyshev polynomials)
the number of terms in which do not dependent on the magnitude of $q$.
In \cite{BremerKummer}, a highly effective numerical method
for constructing  a nonoscillatory solution of  Kummer's equation is described.

The scheme of this paper proceeds by first applying
the algorithm of \cite{BremerKummer} in order to obtain
a nonoscillatory phase function $\alpha$ for (\ref{introduction:second_order}).
We then compute the inverse function $\alpha^{-1}$
(since $\alpha$ is an increasing function, it is invertible).
The roots of a solution $y$ of (\ref{introduction:second_order})
can be easily computing  using $\alpha^{-1}$.
In the event that $y$ is one of the classical orthogonal polynomials, the
values of $y'$ at the roots of $y$
can be calculated and used to construct the weights of the corresponding Gaussian quadrature rule.
The functions  $\alpha$ and $\alpha^{-1}$ are represented by piecewise Chebyshev
expansions whose number of terms is independent of the magnitude of $q$.
Both the time required to compute these expansions and the time required
to evaluate them is independent of the magnitude of $q$.
Moreover, once the computation of $\alpha$ and $\alpha^{-1}$ is completed,
the calculation of each root $t_k$ can be conducted independently. 
In most circumstances, the cost of computing the phase function $\alpha$
and its inverse $\alpha^{-1}$ is small compared to that of calculating
the desired roots of the special function, with the consequence
that the algorithm of this paper admits an effective parallel implementation
(see, for instance, the experiments presented in 
Sections~\ref{section:experiments:jacobi} and \ref{section:experiments:bessel} of this paper).

The remainder of this paper is organized as follows.  In Section~\ref{section:phase},
we briefly review the properties of phase functions. 
Section~\ref{section:nonoscillatory} discusses nonoscillatory phase functions.
In Section~\ref{section:method}, we recount the method of
\cite{BremerKummer} for the numerical computation of nonoscillatory phase
functions.
Section~\ref{section:algorithm} describes an algorithm for the numerical computation
of the roots of special functions satisfying differential equations 
of the form (\ref{introduction:second_order}).
In Section~\ref{section:experiments}, we describe several
numerical experiments conducted to illustrate the properties of the algorithm
of Section~\ref{section:algorithm}.    We conclude with a few remarks
in Section~\ref{section:conclusion}.

\begin{section}{Phase functions and Kummer's equation}

By differentiating (\ref{introduction:u}) twice and adding $q(t)u(t)$ to both
sides of the resulting equation, we see that 
\begin{equation}
u''(t) + q(t) u(t)= 
u(t) 
\left(
 q(t) 
- (\alpha'(t))^2
- \frac{1}{2}\left(\frac{\alpha'''(t)}{\alpha'(t)}\right)
+ \frac{3}{4}
\left(\frac{\alpha''(t)}{\alpha'(t)}\right)^2
\right)
\label{phase:ueq1}
\end{equation}
for all $a < t < b$.  Applying an analogous sequence of steps
to (\ref{introduction:v}) shows that
\begin{equation}
v''(t) + q(t) v(t) = 
v(t) 
\left(
 q(t) 
- (\alpha'(t))^2
- \frac{1}{2}\left(\frac{\alpha'''(t)}{\alpha'(t)}\right)
+ \frac{3}{4}
\left(\frac{\alpha''(t)}{\alpha'(t)}\right)^2
\right)
\label{phase:veq1}
\end{equation}
for all $a < t < b$.  Since $\{u,v\}$ is a basis in the space
of solutions of (\ref{introduction:second_order}), $u$ and $v$
do not simultaneously vanish.
Hence, (\ref{phase:ueq1}) and (\ref{phase:veq1}) imply that
 $\alpha$ is a phase function
for (\ref{introduction:second_order}) if and only if 
its derivative satisfies the second order nonlinear differential equation
\begin{equation}
 q(t) 
- (\alpha'(t))^2
- \frac{1}{2}\left(\frac{\alpha'''(t)}{\alpha'(t)}\right)
+ \frac{3}{4}
\left(\frac{\alpha''(t)}{\alpha'(t)}\right)^2 = 0
\label{phase:kummer}
\end{equation}
on the interval $(a,b)$.  We will refer to (\ref{phase:kummer}) as Kummer's equation, 
after E. E. Kummer who studied it in the 1840s \cite{Kummer}.

If $\alpha$ is a phase function for (\ref{introduction:second_order}), then
any solution $y$ of (\ref{introduction:second_order}) admits a representation
of the form
\begin{equation}
y(t) = 
c_1 \frac{\cos(\alpha(t))}{\sqrt{\alpha'(t)}} +
c_2 \frac{\sin(\alpha(t))}{\sqrt{\alpha'(t)}}.
\label{phase:solution}
\end{equation}

When the coefficient $q$ is of large magnitude, 
the magnitude of $\alpha$ must be large 
(this is a consequence of the Sturm comparison theorem).
In this event, the evaluation of the expression (\ref{phase:solution})
involves the calculation of trigonometric functions of large arguments,
and there is an inevitable loss of precision when such calculations are performed
in finite precision arithmetic.  Nonetheless, acceptable accuracy is obtained
in many cases.  For instance, Section~5.3 of \cite{BremerKummer}
describes an  experiment in which the Bessel function of the first kind of order $10^8$
was evaluated at a large collection of points on the real axis with approximately
ten digits of accuracy.

Moreover, given a phase function $\alpha$ for (\ref{introduction:second_order}),
the roots
\begin{equation}
t_1,\ldots,t_n
\label{phase:roots}
\end{equation}
of a solution $y$ of  (\ref{introduction:second_order}) and the values
\begin{equation}
y'(t_1),\ldots,y'(t_n)
\label{phase:derivative_values}
\end{equation}
of the derivative of $y$ at the roots (\ref{phase:roots}) can be computed without evaluating
trigonometric functions of large arguments and 
the concomitant loss of precision.  To see this, we suppose that $\alpha$ is a phase function
for (\ref{introduction:second_order}) such that $\alpha(a)=0$ and that the values of $y(a)$ 
and $y'(a)$ are known.  An elementary calculation shows that
 the constants $c_1$ and $c_2$ in (\ref{phase:solution})
are given by
\begin{equation}
c_1 = y(a) \sqrt{\alpha'(a)}
\label{phase:c1}
\end{equation}
and
\begin{equation}
c_2 = y(a) \frac{\alpha''(a)}{2 \left(\alpha'(a)\right)^{3/2}} + y'(a) \frac{1}{ \sqrt{\alpha'(a)}}.
\label{phase:c2}
\end{equation}
We let $d_1, d_2$ be the unique pair of real numbers such that  $0 < d_2 \leq \pi$, 
\begin{equation}
c_1 = d_1 \sin(d_2)
\label{phase:d1}
\end{equation}
and
\begin{equation}
c_2 = d_1 \cos(d_2).
\label{phase:d2}
\end{equation}
Then
\begin{equation}
\begin{aligned}
y(t) 
= d_1 \frac{\sin(\alpha(t)+d_2)}{\sqrt{\alpha'(t)}}
\end{aligned}
\label{phase:solution2}
\end{equation}
for all $a \leq t \leq b$.
The expression (\ref{phase:solution2}) is highly conducive to computing the roots of $y$ 
as well as the values of the derivatives of $y$ at the roots of $y$.
Indeed, from (\ref{phase:solution2}) it follows that
\begin{equation}
t_k =  \alpha^{-1} \left( k \pi - d_2 \right) 
\label{phase:inverse_formula}
\end{equation}
for each $k=1,2,\ldots,n$, where $\alpha^{-1}$ denotes the inverse of the
(monotonically increasing function) $\alpha$.
By differentiating (\ref{phase:solution2}), we obtain 
\begin{equation}
y'(t) = d_1 \cos\left(\alpha(t) +d _2\right) \sqrt{\alpha'(t)}
- d_1 \sin\left(\alpha(t) + d_2 \right) \frac{\alpha''(t)}{\left(\alpha'(t)\right)^{3/2}}.
\label{phase:derivative}
\end{equation}
Since
\begin{equation}
\sin(\alpha(t_k)+d_2) = \sin(k\pi) = 0
\end{equation}
and 
\begin{equation}
\cos(\alpha(t_k)+d_2) = \cos(k\pi) = (-1)^k
\end{equation}
for each $k=1,2,\ldots,n$, we see that
\begin{equation}
y'\left(t_k\right) = (-1)^k d_1  \sqrt{\alpha'\left(t_k\right)}
\label{phase:derivative_formula}
\end{equation}
for each $k=1,2,\ldots,n$.  Neither the computation 
of the constants $d_1$ and $d_2$ 
nor the evaluation of formulas
(\ref{phase:inverse_formula}), (\ref{phase:derivative_formula})  
requires the calculation of trigonometric functions of large arguments.

\vskip 1em
\begin{remark}
An obvious modification of the  procedure just described applies
in the case where the values of a solution $y$ of (\ref{introduction:second_order}) and
its derivative are known at arbitrary point in the interval $[a,b]$.
\end{remark}


\label{section:phase}
\end{section}

\begin{section}{Nonoscillatory phase functions}

When $q$ is positive and of large magnitude, almost all phase functions for (\ref{introduction:second_order})
are highly oscillatory.  
Indeed, by differentiating the equation
\begin{equation}
\tan(\alpha(t)) = \frac{v(t)}{u(t)},
\end{equation}
which is easily obtained from (\ref{introduction:u}) and (\ref{introduction:v}),
we see that
\begin{equation}
\alpha'(t) = \frac{W}{\left(u(t)\right)^2 + \left(v(t)\right)^2},
\label{nonoscillatory:alphap}
\end{equation}
where $W$ is the (necessarily constant) Wronskian of the basis $\{u,v\}$.
Since $u$ and $v$ oscillate rapidly when $q$ is positive and of large magnitude
(this is a consequence of the Sturm comparison theorem), 
Formula~(\ref{nonoscillatory:alphap}) 
shows that $\alpha'$, and hence $\alpha$,  will be oscillatory unless some fortuitous cancellation
takes place.

In certain cases, a pair $u, v$ for which such cancellation occurs can be obtained
by finding a solution of (\ref{introduction:second_order}) which is
an element of one of the Hardy spaces  (see, for instance, \cite{Garnett} 
or \cite{Koosis} for an introduction
to the theory of Hardy spaces).    Legendre's differential equation 
\begin{equation}
(1-z^2) y''(z) - 2 z y'(z) + \nu (\nu + 1) y(z) = 0
\label{nonoscillatory:legendre}
\end{equation}
provides one such example.    It can be transformed into the normal form
\begin{equation}
\psi''(z) + \left( \frac{1}{(1-z^2)^2} + \frac{\nu (\nu+1)}{1-z^2}\right) \psi(z) = 0
\label{nonoscillatory:legendre_normal}
\end{equation}
by letting
\begin{equation}
\psi(z) = \sqrt{1-z^2}\ y(z).
\end{equation}
According to  Formula~(9) in Section~3.4 of \cite{HTFI} (see also
Formula~8.834.1 in \cite{Gradshteyn}),
the function $f_\nu$ defined for $x \in \mathbb{R}$ via 
\begin{equation}
f_\nu(x) =  \left( Q_\nu(x)  - i \frac{\pi}{2} P_\nu(x) \right) \sqrt{1-x^2},
\label{nonoscillatory:legendresol}
\end{equation}
where $P_\nu$ and $Q_\nu$ are the Legendre functions of the first and second kinds of degree $\nu$,
respectively, 
  is the boundary value of the solution
\begin{equation}
F_\nu(z) = 
(2z)^{-\nu-1} \sqrt{\pi} \frac{\Gamma(\nu+1)}{\Gamma(\nu+3/2)}
\ \Hypergeometric{2}{1}{\frac{\nu}{2}+1,\frac{\nu}{2}+\frac{1}{2}}{\nu+\frac{3}{2}}{\frac{1}{z^2}}
\sqrt{1-z^2}
\label{nonoscillatory:secondkind}
\end{equation}
of (\ref{nonoscillatory:legendre_normal}).
\begin{figure}[t!!]
\begin{center}
\includegraphics[width=.49\textwidth]{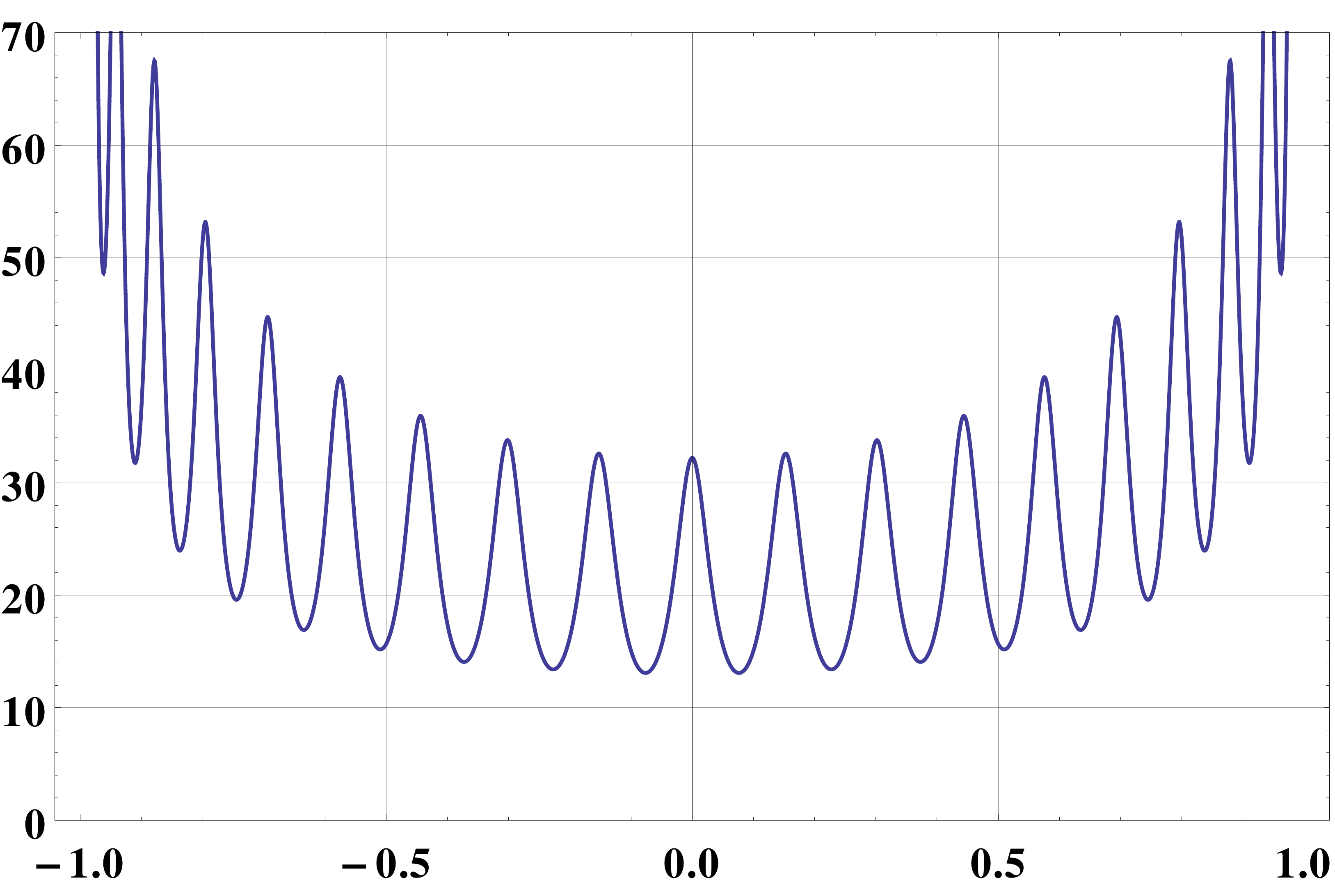}
\hfill
\includegraphics[width=.49\textwidth]{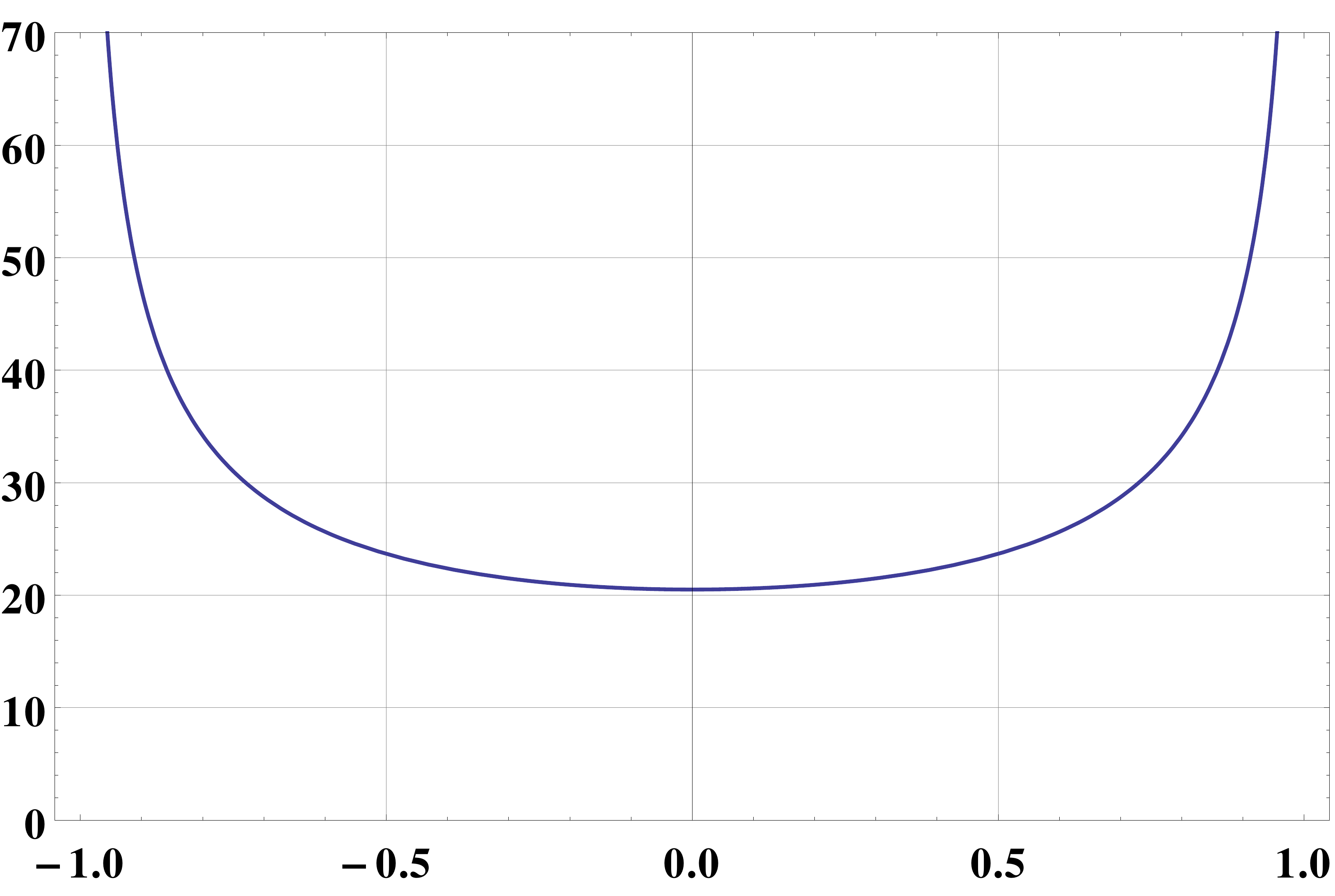}
\end{center}
\caption{\small
On the left is a plot of the derivative of a 
typical (oscillatory) phase function for Legendre's differential equation
when $\nu = 20$.  On the right is a plot of the derivative of a nonoscillatory
phase function for Legendre's differential equation when $\nu=20$.
}
\label{nonoscillatory:figure1}
\end{figure}
The function $F_\nu$ is analytic in the upper half-plane and has no zeros there \cite{Runkel}.
Moreover, although it is not an element of 
one the Hardy spaces on the upper half-plane, its composition with the conformal mapping
\begin{equation}
\tau(z) = i\ \frac{1-z}{1+z}
\end{equation}
of the unit disk onto the upper half-plane is contained in a  Hardy space.
    The imaginary part of the logarithmic derivative
of $f_\nu$ is, of course, the derivative of a phase function for Legendre's differential
equation.  That it is nonoscillatory can be established in a number of ways,
including via the well-known theorem on the factorization of functions in $H^p$  spaces
(which appears as Theorem~5.5 in \cite{Garnett}).  A proof along these
lines will be reported by the author at a later data.
Figure~\ref{nonoscillatory:figure1} depicts
this function 
as well as  the derivative of a typical phase function for (\ref{nonoscillatory:legendre_normal})
 when $\nu = 20$.

This approach can be applied to other equations of interest, including
Bessel's equation (see, for instance,
\cite{Heitman-Bremer-Rokhlin-Vioreanu}), Chebyshev's equation, and the Airy equation.
However, it suffers from at least two significant disadvantages: not every differential
equation of the form (\ref{introduction:second_order}) whose coefficient $q$ is nonoscillatory
can be treated in this fashion
and, perhaps more seriously, even in cases in which it does apply
there is no obvious method for the fast and accurate evaluation of the 
resulting phase functions.

The first of these difficulties is addressed in 
\cite{Heitman-Bremer-Rokhlin} and \cite{Bremer-Rokhlin}.
They contain proofs that, under mild conditions on the 
 coefficient $q$ appearing in (\ref{introduction:second_order}),
there exists
 a nonoscillatory function $\alpha$ such that the functions
(\ref{introduction:u}), (\ref{introduction:v}) approximate
solutions of (\ref{introduction:second_order}) with high accuracy.
The function $\alpha$ is nonoscillatory in the sense that it can be represented
using various series expansions the number of terms of which does not depend
on the magnitude of $q$.

We now state a version of the principal result of 
\cite{Bremer-Rokhlin}
which   pertains to linear ordinary differential equations of the form
\begin{equation}
y''(t) + \lambda^2 q(t) y(t) = 0
\label{nonoscillatory:second_order}
\end{equation}
with $\lambda$ a positive real constant and $q$ a strictly positive function
defined on the real line.  The parameter $\lambda$ is introduced in
order to make rigorous the notion of ``the magnitude of $q$.''  The result 
can be easily applied in cases in which $q$ varies with $\lambda$, as long as
$q$ satisfies the hypotheses of the theorem independent of $\lambda$. 
These requirements are quite innocuous and the theorem can be applied
 to essentially any ordinary differential equation of the form (\ref{introduction:second_order})
with $q$ nonoscillatory.  
See, for instance, the experiments of Section~\ref{section:experiments:artificial}
in which $q$ is taken to be 
\begin{equation}
q(t,\lambda) = 
 \lambda^2\frac{1}{0.1+t^2} + \lambda^{3/2} \frac{\sin(4t)^2}{(0.1+(t-0.5)^2)^4}.
\label{nonoscillatory:q}
\end{equation}
See also Figure~\ref{figure:artificial}, which
 contains a plot of the coefficient defined in Formula~(\ref{nonoscillatory:q}) when $\lambda = 10^5$
as well as a plot of an associated nonoscillatory phase function.

\begin{figure}[b!!]
\begin{center}
\includegraphics[width=.49\textwidth]{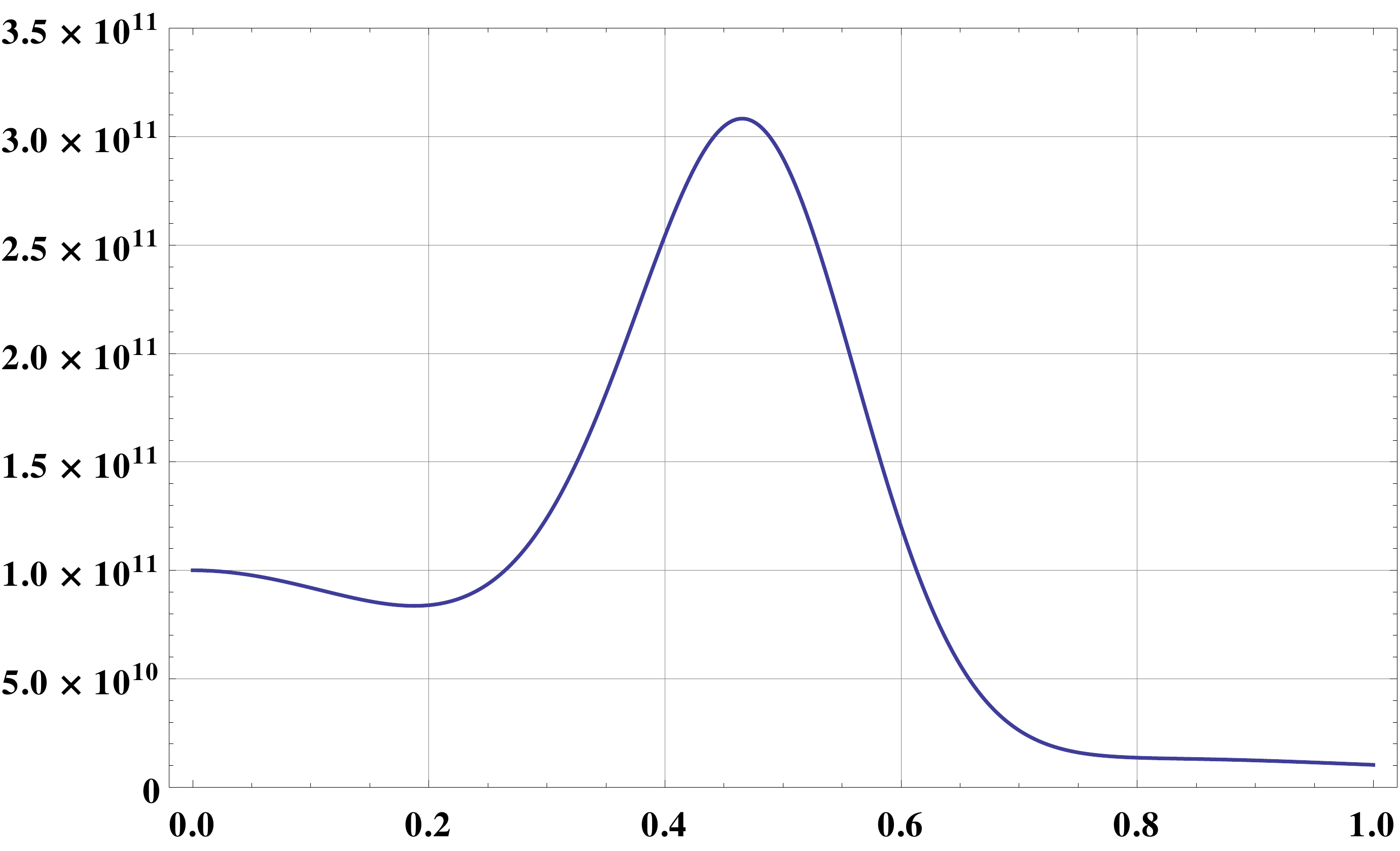}
\hfill
\includegraphics[width=.49\textwidth]{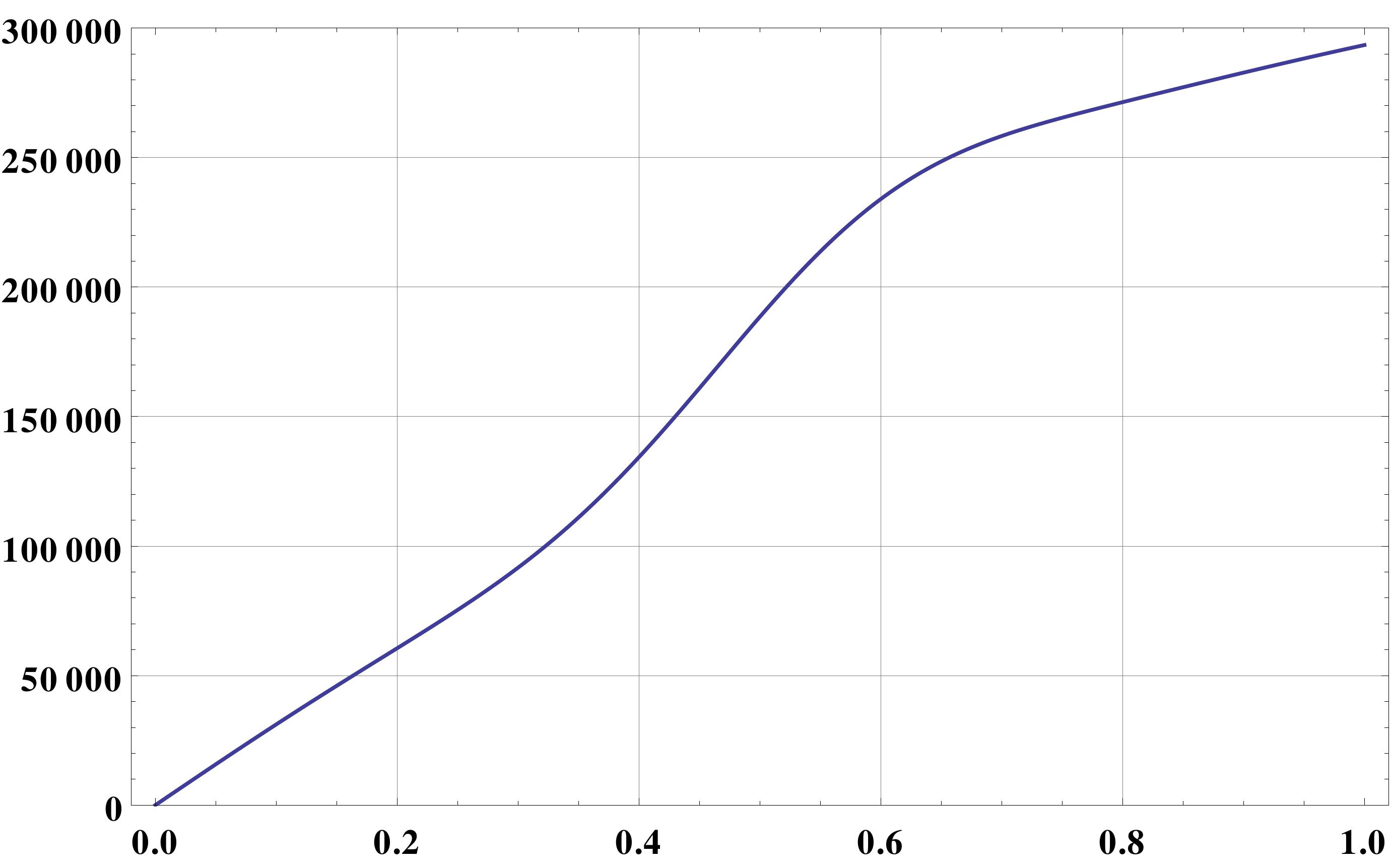}
\end{center}
\caption{\small
On the left is a plot of the function 
(\ref{nonoscillatory:q})
when $\lambda = 10^5$.  On the right is a plot of an associated nonoscillatory
phase function.
}
\label{figure:artificial}
\end{figure}

Similarly, in \cite{Bremer-Rokhlin} it is assumed that
$q$ extends to the real line so that the Fourier transform
can be used to quantify the notion of ``nonoscillatory function.''  
In most cases of interest, the coefficient $q$ has either zeros,
singularities or both at the endpoints of the interval on which
(\ref{nonoscillatory:second_order}) is given.
The coefficient $q$
in  (\ref{nonoscillatory:legendre_normal}), for instance, has poles at the points
$\pm 1$; nonetheless 
it  can be  represented to high relative accuracy via (for instance) local
expansions in polynomials given on ``graded intervals.'' That is,
on a collection of subintervals which become smaller as they approach
the points $\pm 1$.  By extending each local expansion to the real line,
the results of \cite{Bremer-Rokhlin} can be brought to bear on the problem.
Equations with turning points (i.e., locations where the coefficient $q$ vanishes)
can be dealt with in a similar fashion.

Several definitions are required before stating the principal result of \cite{Bremer-Rokhlin}.
We say that an infinitely differentiable function $\varphi: \mathbb{R} \to \mathbb{R}$
is a Schwartz function if
\begin{equation}
\sup_{t\in\mathbb{R}} |t^i \varphi^{(j)}(t)| < \infty
\end{equation}
for all pairs $i,j$ of nonnegative integers, and we denote the set  of all Schwartz functions 
 by  $\Sr$.  We define functions  $\tilde{p}$ and $x$  via the formulas
\begin{equation}
\tilde{p}(t) = \frac{1}{q(t)} \left( \frac{5}{4} \left(\frac{q'(t)}{q(t)}\right)^2 -
\frac{q''(t)}{q(t)} \right)
\label{nonoscillatory:definition_of_ptilde}
\end{equation}
and
\begin{equation}
x(t) = \int_a^t \sqrt{q(u)}\ du.
\label{nonoscillatory:definition_of_x}
\end{equation}
Since $q$ is strictly positive, $x(t)$  is monotonically increasing and hence invertible.
We define the function $p$ via the formula
\begin{equation}
p(x) = \tilde{p}(t(x));
\label{nonoscillatory:definition_of_p}
\end{equation}
that is, $p$ is the composition of $\tilde{p}$ with the inverse of the
function $x$ defined via (\ref{nonoscillatory:definition_of_x}). 
The relationship between the function $p$ and $q$ is 
ostensibly complicated, but, in fact,
\begin{equation}
p(x)= 2 \left\{t,x\right\},
\label{nonoscillatory:p_schwarzian}
\end{equation}
where $\{t,x\}$ denotes the Schwarzian derivative of the inverse of the function
$t(x)$ .  That is, $\{f,x\}$ is defined via the formula
\begin{equation}
\{f,x\} = \frac{f'''(x)}{f'(x)} - \frac{3}{2} \left(\frac{f''(x)}{f'(x)}\right)^2.
\end{equation}
A derivation of (\ref{nonoscillatory:p_schwarzian}) can be found in 
Section~3 of \cite{Bremer-Rokhlin}.
The following theorem is a consequence of Theorem~12 in \cite{Bremer-Rokhlin}.

\vskip 1em
\begin{theorem}
Suppose that the function $\tilde{p}$ defined via (\ref{nonoscillatory:definition_of_p})
is an element of the Schwartz space $\Sr$, 
 that there exist positive real numbers $\Gamma$ and $\mu$ such that
\begin{equation}
\left|\widehat{p}(\xi)\right| \leq \Gamma \exp(-\mu |\xi|)
\ \ \mbox{for all} \ \xi\in\mathbb{R},
\label{nonoscillatory:phat}
\end{equation}
and that  $\lambda$ is a positive real number such that
\begin{equation}
\lambda > 2 \max\left\{\frac{1}{\mu},\Gamma\right\}.
\label{nonoscillatory:lambda}
\end{equation}
Then there exist functions  $\nu$ and $\delta$ in $\Sr$
such that
\begin{equation}
\left|\widehat{\delta}(\xi)\right| \leq 
\frac{\Gamma}{\lambda^2} \left( 1+ \frac{2\Gamma}{\lambda}\right) 
\exp\left(-\mu\left|\xi\right|\right)
\ \ \ \mbox{for all} \ \ \xi \in \mathbb{R},
\label{nonoscillatory:deltahat}
\end{equation}
\begin{equation}
\|\nu\|_\infty \leq \frac{\Gamma}{2\mu} \left(1 + \frac{4\Gamma}{\lambda}\right)
\exp\left(-\mu\lambda\right),
\label{nonoscillatory:nubound}
\end{equation}
and the function $\alpha$ defined via the formula
\begin{equation}
\alpha(t) = \lambda \sqrt{q(t)} \int_a^t   \exp\left(\frac{1}{2} \delta(u) \right)\ du
\label{nonoscillatory:alpha_definition}
\end{equation}
satisfies the nonlinear differential equation
\begin{equation}
\left(\alpha'(t)\right)^2 = \lambda^2 
\left(\frac{\nu(t)}{4\lambda^2}+1\right)q(t) - \frac{1}{2}\frac{\alpha'''(t)}{\alpha'(t)}
+ \frac{3}{4} \left(\frac{\alpha''(t)}{\alpha'(t)}\right)^2.
\label{nonoscillatory:alpha2}
\end{equation}
\label{main_theorem}
\end{theorem}

It follows from (\ref{nonoscillatory:alpha2}) that
 $\alpha$ is a phase function for the perturbed second order linear ordinary
differential equation
\begin{equation}
y''(t) + \lambda^2 \left(1+\frac{\nu(t)}{4\lambda^2}\right)q(t)y(t) = 0,
\label{nonoscillatory:second_order2}
\end{equation}
and we conclude from this fact and (\ref{nonoscillatory:nubound}) 
that when $\alpha$ is inserted into formulas 
(\ref{introduction:u}) and (\ref{introduction:v}), 
the resulting functions approximate solutions of 
(\ref{introduction:second_order}) with accuracy on the order of 
\begin{equation}
\frac{1}{\mu\lambda}\exp(-\mu \lambda).
\end{equation}
From  (\ref{nonoscillatory:deltahat})  we see that
the Fourier transform of $\delta$ is bounded by $\exp(-\mu |\xi|)$
for sufficiently large $\lambda$.  Among other things, this implies
that the function $\delta$ can be represented on the interval $[a,b]$
on which (\ref{introduction:second_order}) is given
via various series expansions (e.g., as a Chebyshev or Legendre expansion)
using a number of terms which does not depend on $\lambda$.  
Plainly, the function $\alpha$ defined via
(\ref{nonoscillatory:alpha_definition})
has this property as well.  It is in this sense that 
 $\alpha$ is nonoscillatory ---
it can be accurately represented
using various series expansions
whose number of terms does not depend on the parameter $\lambda$ (which is
a proxy for the magnitude of the coefficient $q$).
Theorem~\ref{main_theorem} can be summarized by saying that
nonoscillatory phase functions 
represent the  solutions of (\ref{nonoscillatory:second_order}) with 
 $\mathcal{O}\left((\mu\lambda)^{-1}\exp(-\mu\lambda)\right)$ accuracy
using  $\mathcal{O}(1)$ terms.  
Using a nonoscillatory phase function,
Formulas~(\ref{phase:inverse_formula})
and (\ref{phase:derivative_formula}) can be evaluated in $O(1)$ operations. 
This is in contrast to the $O(\lambda)$ operations required
to evaluate them when $\alpha$ is a typically, oscillatory phase function
for (\ref{nonoscillatory:second_order}).

By a slight abuse of terminology, throughout the rest of this article
we will refer to the 
 the function $\alpha$ defined via (\ref{nonoscillatory:alpha_definition})
as the nonoscillatory solution of Kummer's equation
and as the nonoscillatory phase function for
Equation~(\ref{nonoscillatory:second_order}).
Note that this function is neither an exact solution
of (\ref{phase:kummer})
(although it approximates
a solution of that equations with error which decays exponentially with
$\lambda$), nor is it unique (but it is 
the  one and only one 
nonoscillatory phase function associated with Theorem~\ref{main_theorem}).

\label{section:nonoscillatory}
\end{section}

\begin{section}{A practical method for the computation of nonoscillatory phase functions}

The method used in \cite{Heitman-Bremer-Rokhlin} and \cite{Bremer-Rokhlin} to prove
the existence of nonoscillatory phase function is constructive and could serve
as the basis of a numerical algorithm for their computation.  
However, such an approach would require that $q$ be explicitly extended
to the real line as well as knowledge of the first two derivatives of $q$.
We now describe a minor variant of the algorithm of \cite{BremerKummer}.
Under the hypotheses of Theorem~\ref{main_theorem}, it results in
a nonoscillatory function $\alpha$ which represents
solutions of (\ref{nonoscillatory:second_order})
with accuracy on the order of  $\exp\left(-\frac{1}{2}\mu\lambda\right)$,
where $\mu$ is the constant appearing in Theorem~\ref{main_theorem}.
The method of \cite{BremerKummer} 
has the advantage that it only requires knowledge of 
the coefficient $q$ on the interval $[a,b]$.  Note that, as is the case with Theorem~\ref{main_theorem},
there is no difficulty in treating equations in which the coefficient $q$ varies with $\lambda$,
assuming that $q$ satisfies the hypotheses of Theorem~\ref{main_theorem} independent of $\lambda$.

\begin{figure}[b!!]
\begin{center}
\includegraphics[width=.46\textwidth]{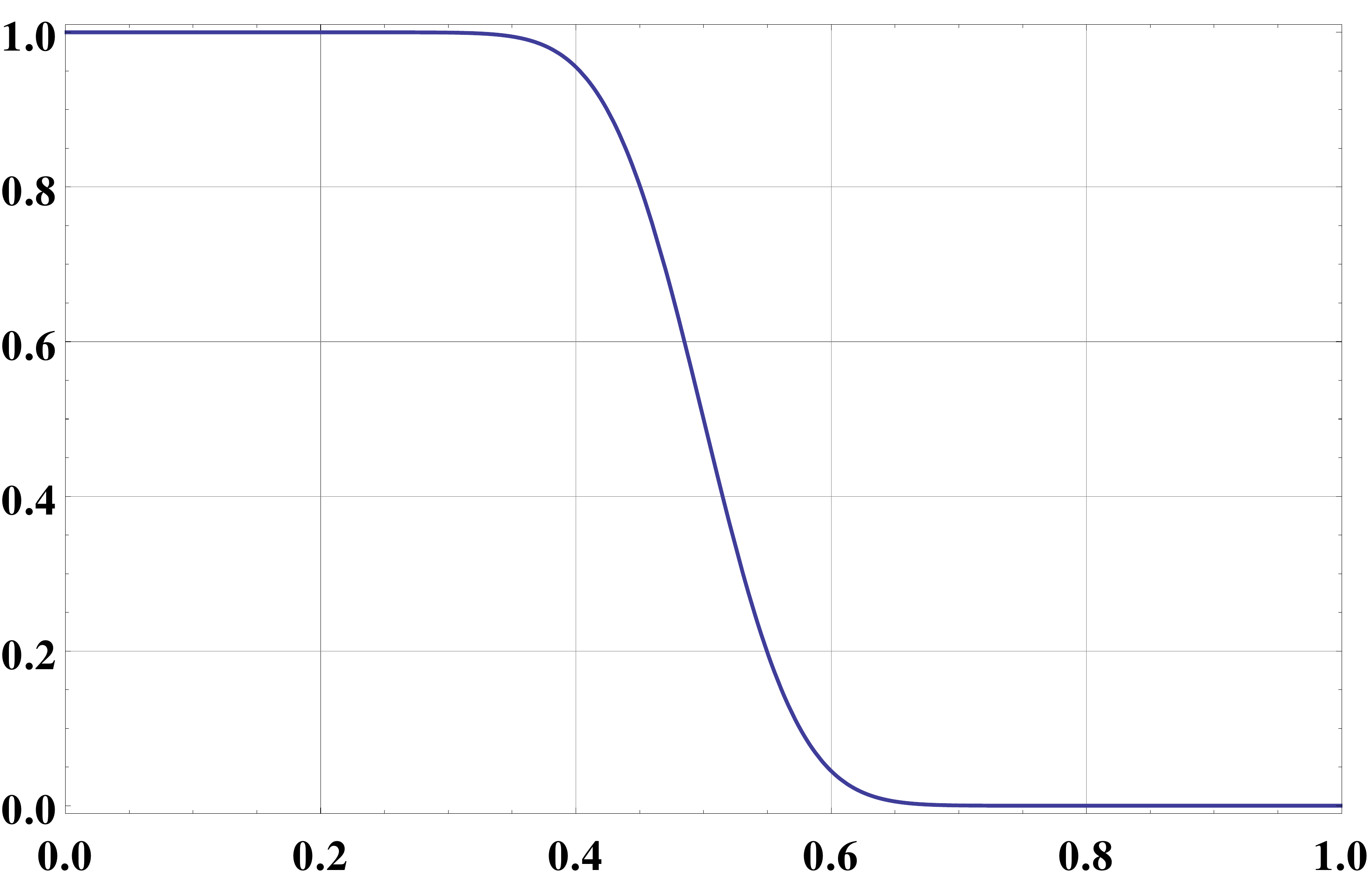}
\hfill
\includegraphics[width=.49\textwidth]{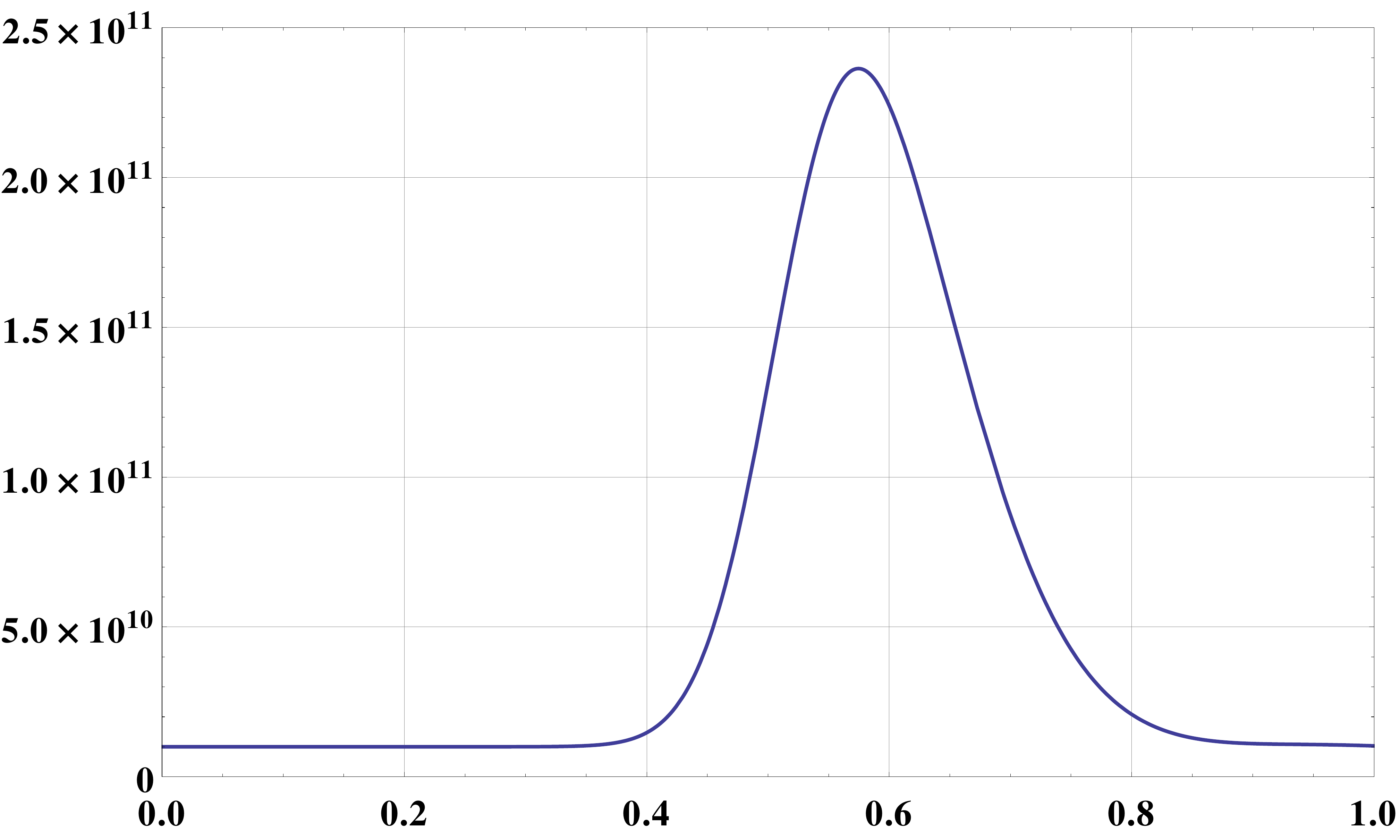}
\caption{\small
On the left is the function $\phi$ defined by (\ref{method:phi})
when $a=0$ and $b=1$.
On the right is the windowed version of the coefficient $q$ 
given by (\ref{nonoscillatory:q}) with $\lambda=10^5$.  
It agrees with the original coefficient $q$, 
a plot of which 
appears in Figure~\ref{figure:artificial},
on the rightmost quarter of the interval $[0,1]$ and is constant
on the leftmost quarter of $[0,1]$.
}
\label{method:figure1}
\end{center}
\end{figure}

The algorithm  proceeds as follows.
First, a windowed version $\tilde{q}$ of
$q$ which closely agrees with $q$ on the rightmost quarter of the interval $[a,b]$
and is approximately equal to the constant $1$ on the leftmost
quarter of $[a,b]$ is constructed.  More specifically, $\tilde{q}$ is defined via
the formula
\begin{equation}
\tilde{q}(t) = \phi(t) + (1-\phi(t)) q(t),
\end{equation}
where $\phi$ is given by
\begin{equation}
\phi(t) = \frac{1}{2} 
\left( 
\erf\left(\frac{24}{b-a} \left(t+\frac{a+b}{2}\right)\right) - 
\erf\left(\frac{24}{b-a} \left(t-\frac{a+b}{2}\right)\right) 
\right).
\label{method:phi}
\end{equation}
We  use the analytic windowing function $\phi$ to construct $\tilde{q}$ 
(as opposed to a windowing function which is infinitely differentiable and compactly supported)
so that we can apply Theorem~\ref{main_theorem} to the windowed
version $\tilde{q}$ of $q$ (its hypotheses imply that $q$ is analytic
in a strip containing the real line).
The constants in (\ref{method:phi}) are set so that
\begin{equation}
1 - \phi(t) < 10^{-16} \ \ \mbox{for all} \ \ t \leq \frac{3a+b}{4}
\end{equation}
and 
\begin{equation}
 \phi(t) < 10^{-16} \ \ \mbox{for all} \ \ t \geq \frac{a+3b}{4}.
\end{equation}
The left side of Figure~\ref{method:figure1} shows a plot of the function $\phi$
when $a=0$ and $b=1$.  Next, a solution $\tilde{\alpha}'$ of
 the initial value problem
\begin{equation}
\left\{
\begin{aligned}
 \lambda^2 \tilde{q}(t) - \left(\tilde{\alpha}'(t)\right)^2 - \frac{1}{2}\frac{\tilde{\alpha}'''(t)}
{\tilde{\alpha}'(t)}
+ \frac{3}{4}\left(\frac{\tilde{\alpha}''(t)}{\tilde{\alpha}'(t)}\right)^2 &= 0 
\ \ \mbox{for all} \ \ a \leq t \leq b \\
\tilde{\alpha}'(a) &= \lambda \\
\tilde{\alpha}''(a) & = 0
\end{aligned}
\right.
\label{method:bvp1}
\end{equation}
is obtained.  The initial conditions $\tilde{\alpha}'(a) = \lambda$ and 
$\tilde{\alpha}''(a)=0$
in  (\ref{method:bvp1}) are chosen because the nonoscillatory
phase function for  the equation
\begin{equation}
y''(t) + \lambda^2 \tilde{q}(t) y(t) = 0
\label{method:perturbed}
\end{equation}
whose existence is ensured by Theorem~\ref{main_theorem}
behaves as 
\begin{equation}
\lambda t + O\left(\exp\left(-\frac{1}{2}\mu\lambda\right)\right)
\end{equation}
on the leftmost quarter of the interval $[a,b]$ where
$\lambda^2 \tilde{q}$ is approximately equal to the constant $\lambda^2$.
This estimate is proven in  \cite{BremerKummer}.  It is also shown in
\cite{BremerKummer} that on the rightmost quarter of the interval [$a,b]$
the difference between the nonoscillatory phase function 
for Equation~(\ref{nonoscillatory:second_order}) 
and the solution $\tilde{\alpha}$ 
of (\ref{method:bvp1}) is on the order of 
$\exp\left(-\frac{1}{2} \mu \lambda\right)$
as is the difference between $\tilde{\alpha}'$
and the derivative of the nonoscillatory phase function for 
Equation~(\ref{nonoscillatory:second_order}).
It follows that by solving the initial value problem
(\ref{method:bvp1}) we approximate
the values of the nonoscillatory phase function for (\ref{nonoscillatory:second_order})
and its derivative at the right-hand endpoint $b$ 
with accuracy on the order of $\exp(-\frac{1}{2} \mu \lambda)$.
  The algorithm concludes by solving the terminal value
problem
\begin{equation}
\left\{
\begin{aligned}
 \lambda^2 q(t) - \left(\alpha'(t)\right)^2 - \frac{1}{2}\frac{\alpha'''(t)}{\alpha'(t)}
+ \frac{3}{4}\left(\frac{\alpha''(t)}{\alpha'(t)}\right)^2 &= 0 
\ \ \mbox{for all} \ \ a \leq t \leq b\\
\alpha'(b) &= \tilde{\alpha}'(b) \\
\alpha''(b) &= \tilde{\alpha}''(b). \\
\end{aligned}
\right.
\label{method:bvp2}
\end{equation}
Standard results on the continuity of solutions of ordinary differential
equations (see, for instance, \cite{Coddington-Levinson}) together
with the preceding estimate
imply that  the solution $\alpha$ of (\ref{method:bvp2}) approximates
the nonoscillatory phase function
for Equation~(\ref{nonoscillatory:second_order}) on the interval $[a,b]$
with accuracy on the order of
$\exp\left(-\frac{1}{2} \mu \lambda\right)$.

\label{section:method}
\end{section}

\begin{section}{An algorithm for the computation of the roots of special functions}

In this section, we  describe an algorithm for computing the roots
\begin{equation}
t_1< \ t_2 < \ldots < t_n
\label{algorithm:roots}
\end{equation}
on the interval $[a,b]$ of a solution $y$ of Equation~(\ref{nonoscillatory:second_order}), 
as well as the values of the derivative of $y$ at the points  (\ref{algorithm:roots}).
It takes as inputs the value of $\lambda$, a subroutine for evaluating the coefficient
$q$ in (\ref{introduction:second_order}), the values of the function
$y$ and its derivative at the left-hand end point $a$ of the interval
$[a,b]$ on which (\ref{introduction:second_order}) is given, a positive
integer $k$, and a partition
\begin{equation}
a = \gamma_1 < \gamma_2 < \ldots < \gamma_m < \gamma_{m+1} = b
\label{algorithm:partition}
\end{equation}
of the interval $[a,b]$.    As with the analyses of Section~\ref{section:nonoscillatory}
and \ref{section:method}, the algorithm described here can be easily adapted
to the case in which the coefficient $q$ varies with $\lambda$.

The algorithm proceeds as follows:
\begin{enumerate}

\item
We construct a nonoscillatory solution $\tilde{\alpha}$ of the initial value problem
(\ref{method:bvp1}).  Then, having obtained the values of $\tilde{\alpha}(b)$ and
$\tilde{\alpha}(b)$, we solve the terminal value problem (\ref{method:bvp2}).  
The resulting function $\alpha$ is 
a nonoscillatory solution of Kummer's equation.

\item
We next construct the inverse function $\alpha^{-1}$ of the monotonically increasing
function $\alpha$ 
on the interval $[\alpha(a), \alpha(b)]$ via Newton's method.

\item
Then, we calculate the coefficients $c_1$ and $c_2$ such that
\begin{equation}
y(t) = 
c_1 \frac{\cos(\alpha(t))}{\sqrt{\alpha'(t)}} +
c_2 \frac{\sin(\alpha(t))}{\sqrt{\alpha'(t)}}
\end{equation}
using Formulas~(\ref{phase:c1}) and (\ref{phase:c2}).  We determine $d_1$ and $d_2$ such that
\begin{equation}
y(t) = d_1 \frac{\sin(\alpha(t)+d_2)}{\sqrt{\alpha'(t)}}
\end{equation}
using the relations (\ref{phase:d1}), (\ref{phase:d2}).

\item
Finally, we use Formula~(\ref{phase:inverse_formula}) to calculate the $k^{th}$ root $t_k$
of y on the interval $[a,b]$ and (\ref{phase:derivative_formula}) to calculate the
value of $y'$ at $t_k$
\end{enumerate}
Of course, once Steps~1-4 have been completed, the roots $t_k$ 
and corresponding values $y'(t_k)$ can be  computed in any order
using Formulas~(\ref{phase:inverse_formula}) and (\ref{phase:derivative_formula}).

The function $\alpha$ is represented by its values
at the $k$-point Chebyshev grids on each of the $m$  subintervals 
\begin{equation}
\left[\gamma_1,\gamma_{2}\right],\ \left[\gamma_2,\gamma_{3}\right],\ \ldots, \
\left[\gamma_m,\gamma_{m+1}\right].
\label{algorithm:subintervals}
\end{equation}
The $k$-point Chebyshev grid on $[\gamma_i,\gamma_{i+1}]$ is the set
\begin{equation}
\left\{
\frac{\gamma_{i+1}+\gamma_i}{2} + 
\frac{\gamma_{i+1}- \gamma_i}{2} 
\cdot \cos\left(\frac{j \pi}{k-1}\right) :\ j =0,1,\ldots,k-1\right\}.
\label{algorithm:chebyhev_grid}
\end{equation}
Given the values of a polynomial of degree $k-1$ at the points (\ref{algorithm:chebyhev_grid}),
its value at any point on the
interval $[\gamma_i,\gamma_{i+1}]$ can be calculated  in a numerically stable fashion
in $\mathcal{O}(k)$ operations  via barycentric interpolation.
 Moreover, assuming the partition
(\ref{algorithm:partition}) is properly chosen, 
 barycentric interpolation can be used to approximate value of $\alpha$ at
any  point $x$ in $[a,b]$ given its values at the Chebyshev nodes 
on the subintervals (\ref{algorithm:subintervals}).  This can be done in
 $\mathcal{O}\left(k + \log_2(m)\right)$ operations by first finding the subinterval
containing $x$ using a binary search and then performing barycentric interpolation
on that subinterval.  See, for instance, \cite{Trefethen}
for an extensive discussion of Chebyshev polynomials and interpolation.
In the experiments described in Section~\ref{section:experiments} of this paper, $k$
took on various values between $5$ and $30$.

The problems (\ref{method:bvp1}) and (\ref{method:bvp2}) are stiff when 
$\lambda$ is large and
an appropriately chosen method must be used to solve them numerically.
We compute the values of the solutions at the $k$-point Chebyshev grid
on each of the subintervals $[\gamma_i,\gamma_{i+1}]$
by first constructing a low-accuracy approximation to the solution 
via the implicit trapezoidal method
(see, for instance, \cite{Iserles}), and then applying the Newton-Kantorovich method
(which is discussed in Chapter~5 of \cite{ZeidlerI}, among many other sources).
We now describe the specifics of the Newton-Kantorovich method in more detail.
By multiplying both sides  of Kummer's equation
\begin{equation}
 \lambda^2 q(t) - \left(\alpha'(t)\right)^2 - \frac{1}{2}\frac{\alpha'''(t)}{\alpha'(t)}
+ \frac{3}{4}\left(\frac{\alpha''(t)}{\alpha'(t)}\right)^2 = 0 
\label{algorithm:kummer}
\end{equation}
by $\alpha'$, letting $\beta(t) = \alpha'(t)$ and rearranging the terms of the
resulting equation, we arrive at 
\begin{equation}
\beta''(t) + 2(\beta(t))^3- 2 \lambda^2 q(t) \beta(t)  - \frac{3}{2} \frac{\beta'(t)}{\beta(t)} = 0.
\label{algorithm:beta}
\end{equation}
In each iteration of Newton-Kantorovich method, the given approximation
$\beta_0$ of the solution of (\ref{algorithm:beta}) is  updated by solving
the linearized equation
\begin{equation}
\delta''(t) +  
\left(3 \frac{\beta_0'(t) }{\beta(t)}\right) \delta'(t)  + 
\left(6 (\beta_0(t))^2 + \frac{3}{2} \frac{(\beta'(t))^2}{\left(\beta(t)\right)^2}
- 2 q(t) \right) \delta(t) = r(t),
\label{algorithm:linearized}
\end{equation}
where
\begin{equation}
r(t) = 
-\beta_0''(t) - 2(\beta_0(t))^3  
+ 2 \lambda^2 q(t) \beta_0(t)  + \frac{3}{2} \frac{\beta_0'(t)}{\beta_0(t)},
\end{equation}
for $\delta$ and taking the new approximation of the solution of (\ref{algorithm:beta})
to be $\beta_0 + \delta$.
Equation~(\ref{algorithm:linearized})
is solved via a variant of  the spectral method of \cite{GreengardTwoPoint}.  Newton-Kantorovich
iterations are continued until no further improvement in the solution $\beta$ is 
obtained.  A solution $\alpha$ of Kummer's equation (\ref{algorithm:kummer}
 is obtained from $\beta$ through the formula
\begin{equation}
\alpha(t) = \int_a^t \beta(u)\ du.
\end{equation}
%


The inverse phase function $\alpha^{-1}$ is represented via its values
on the $k$-point Chebyshev grids on each of the $m$ subintervals
\begin{equation}
\left[ \alpha(\gamma_1), \alpha(\gamma_2) \right],\
\left[ \alpha(\gamma_2), \alpha(\gamma_3) \right],\  
\ldots\ 
\left[ \alpha(\gamma_m), \alpha(\gamma_{m+1}) \right].
\label{algorithm:inverse_subintervals}
\end{equation}
There are $m(k-1)+ 1$ such points (since the set (\ref{algorithm:chebyhev_grid}) includes
the endpoints of each interval), and we denote them by
\begin{equation}
\rho_1 < \rho_2 < \ldots < \rho_{m(k-1)+1}.
\label{algorithm:points}
\end{equation}
The values of $\alpha^{-1}$ at the endpoints of the intervals 
(\ref{algorithm:inverse_subintervals})
are (obviously) known.  We calculate its values at the remaining points via Newton's method.
Since $\alpha$ is monotonically increasing, we start the process
by evaluating $\alpha^{-1}$ at $\rho_{m(k-1)}$, then
at $\rho_{m(k-1)-1}$, and so on.  As in the case of $\alpha$,  
once the values of $\alpha^{-1}$ at the points (\ref{algorithm:points})
are known, the value of $\alpha^{-1}$ at any point on the interval
$[\alpha(a), \alpha(b)]$
can be approximated in $\mathcal{O}\left(k+ \log_2(m)\right)$ operations 
via barycentric interpolation.   

In some cases, an appropriate collection of subintervals (\ref{algorithm:subintervals})
is not known {\it a priori} and must be  determined through an adaptive procedure.
In that event, we  construct an initial set of subintervals  by adaptively discretizing
the function $\sqrt{q(t)}$.  We use $\sqrt{q(t)}$ as a starting point 
for the representation of $\alpha'$ because of the classical estimate
\begin{equation}
\alpha'(t) = \lambda \sqrt{q(t)} + O \left(1\right);
\label{algorithm:wkb}
\end{equation}
see, for instance, \cite{Fedoryuk}.  We note too that
when inserted into (\ref{introduction:u}) and (\ref{introduction:v}), the
crude approximation
\begin{equation}
\alpha(t) \approx \lambda \int_a^t \sqrt{q(u)}\ du
\end{equation}
 gives rise to the first order
WKB approximations of the solutions of (\ref{nonoscillatory:second_order}).
The adaptive discretization of $\sqrt{q(t)}$ proceeds by
recursively subdividing the interval $[a,b]$ 
until  each of the resulting subintervals $[a_0,b_0]$
satisfy the following property.  Suppose that
\begin{equation}
\sum_{l=0}^{k-1}  c_l T_l(t)
\label{algorithm:chebyshev_expansion}
\end{equation}
is the Chebyshev expansion of the polynomial of degree $k-1$ which interpolates
 $\sqrt{q(t)}$ at the nodes of the $k$-point Chebyshev grid on the subinterval
$[a_0,b_0]$, and that 
\begin{equation}
c_{\mbox{\tiny max}} = \max\{|c_0|,|c_1|,\ldots,|c_{k-1}|\}.
\end{equation}
We  subdivide $[a_0,b_0]$ if any of the quantities  
\begin{equation}
 \frac{\left|c_{\lceil k/2 \rceil}\right|}{c_{\mbox{\tiny max}}}, \ldots,
 \frac{\left|c_{k-1}\right|}{c_{\mbox{\tiny max}}}
\end{equation}
%
are greater than a specified value (which was taken 
to be $10^{-13}$ in the experiments described in Section~\ref{section:experiments}
of this paper).    

The collection of subintervals obtained by discretizing
$\sqrt{q(t)}$ might not suffice
for representing the solution $\alpha$ of Kummer's equation.
Accordingly, we follow the following procedure while 
solving the problems (\ref{method:bvp1}) and (\ref{method:bvp2}).
After using the Newton-Kantorovich method as described above to approximate
the values of the solution $f$ of one of these problems
on a subinterval $[a_0,b_0]$, we compute the coefficients 
$c_0,c_1,\ldots,c_{k-1}$ in the expansion (\ref{algorithm:chebyshev_expansion})
of the polynomial of degree $k-1$ interpolating $f$ at the nodes of the $k$-point Chebyshev
grid on $[a_0,b_0]$.  We apply the same criterion as before in order to decide
whether to divide the interval $[a_0,b_0]$ or not; that is,
if the relative magnitude of one the trailing 
coefficients in that expansion is too large, then the interval
$[a_0,b_0]$  is split in half and the same procedure is applied, recursively,
to each of the two resulting subintervals of $[a_0,b_0]$.

\vskip 1em
\begin{remark}
It is possible that the collection 
of subintervals (\ref{algorithm:inverse_subintervals}) is insufficient
to represent the function $\alpha^{-1}$; however, we have never seen
this occur in practice.  Indeed, an early version of the algorithm of 
this paper adaptively discretized $\alpha^{-1}$ independently of $\alpha$, 
but this code was removed as it was never necessary.
\end{remark}

\vskip 1em
\begin{remark}
The procedure for the numerical solution of the boundary value 
problems (\ref{method:bvp1}) and (\ref{method:bvp2}) described here
was chosen for its robustness and its ability to solve extremely
stiff ordinary differential equations to high accuracy.  In the 
 numerical experiments of Section~\ref{section:experiments:jacobi},
for instance, values of $\lambda$ as large as $10^{12}$ are considered.
When the value of  $\lambda$ is somewhat  smaller, 
faster methods 
--- such as the spectral deferred correction method of \cite{Dutt-Greengard-Rokhlin}
---  may be used in place of the  technique described here.
\end{remark}


\label{section:algorithm}
\end{section}

\begin{section}{Numerical experiments}

In this section, we describe several numerical experiments which were conducted to
illustrate the performance of the algorithm of this article.
Our code  was written in Fortran 95 using OpenMP extensions   and
compiled  with the Intel Fortran Compiler version 16.0.0.
All calculations were carried out 
on a desktop computer equipped with 28 
Intel Xeon E5-2697  processor cores running at 2.6 GHz and 512 GB of memory.
A version of the code used to conduct these experiments is available 
at the author's website: \url{http://www.math.ucdavis.edu/~bremer/code.html}.

\begin{subsection}{An artificial example}

For various values of $\lambda$, we computed the roots of the solution of the initial value
problem
\begin{equation}
\left\{
\begin{aligned}
y''(t) + q(t,\lambda) y(t) &= 0 \ \ \mbox{for all} \ \ 0 \leq t \leq 1 \\
y(0)   &= 0\\
y'(0)  &= \lambda,
\end{aligned}
\right.
\label{art:eq}
\end{equation}
where $q$ is defined via the formula
\begin{equation}
q(t,\lambda) = 
 \lambda^2\frac{1}{0.1+t^2} + \lambda^{3/2} \frac{\sin(4t)^2}{(0.1+(t-0.5)^2)^4}.
\label{art:q}
\end{equation}

\begin{table}[h!!!!]
\begin{center}
\begin{tabular}{c@{\hspace{1em}}crc}
\toprule
\addlinespace[.5em]
$\lambda$             &  Phase function       & Number of            & Root calculation             \\
                      &  time                 & roots in $[0,1]$     & time                        \\
\midrule
\addlinespace[.75em]
$10^3$    & 7.03\e{-02} & \num[group-separator={,}]{2096} & 2.66\e{-04} \\
\addlinespace[.25em]
$10^4$    & 3.49\e{-02} & \num[group-separator={,}]{13339} & 2.05\e{-03} \\
\addlinespace[.25em]
$10^5$    & 3.18\e{-02} & \num[group-separator={,}]{93398} & 1.47\e{-02} \\
\addlinespace[.25em]
$10^6$    & 2.88\e{-02} & \num[group-separator={,}]{736207} & 9.78\e{-02} \\
\addlinespace[.25em]
$10^7$    & 2.61\e{-02} & \num[group-separator={,}]{6476851} & 7.07\e{-01} \\
\addlinespace[.25em]
$10^8$    & 3.00\e{-02} & \num[group-separator={,}]{61289533} & 5.02\e{+00} \\
\addlinespace[.25em]
$10^9$    & 3.59\e{-02} & \num[group-separator={,}]{600685068} & 4.63\e{+01} \\
\addlinespace[.25em]
\bottomrule
\end{tabular}
\end{center}
\caption{The results of the experiment of Section~\ref{section:experiments:artificial}.
These calculations were performed on a single processor core.  All times are in
seconds.}
\label{art:table}
\end{table}

  The adaptive version of the algorithm of Section~\ref{section:algorithm}
was used and the parameter $k$ was taken to be $16$.
Table~\ref{art:table} presents the results.  There, each row
corresponds to one value of $\lambda$ and reports the time (in seconds) required
to construct the nonoscillatory phase function and its inverse,
the number of roots of the solution of (\ref{art:eq}) in
the interval $[0,1]$, and the time required to compute the roots.
 Figure~\ref{figure:artificial}
displays the coefficient (\ref{art:q}) when $\lambda = 10^5$,
as well as a plot of an associated nonoscillatory phase function.
Figure~\ref{figure:art2} shows the inverse of that nonoscillatory phase function.

\begin{figure}[h!!]
\begin{center}
\includegraphics[width=.6\textwidth]{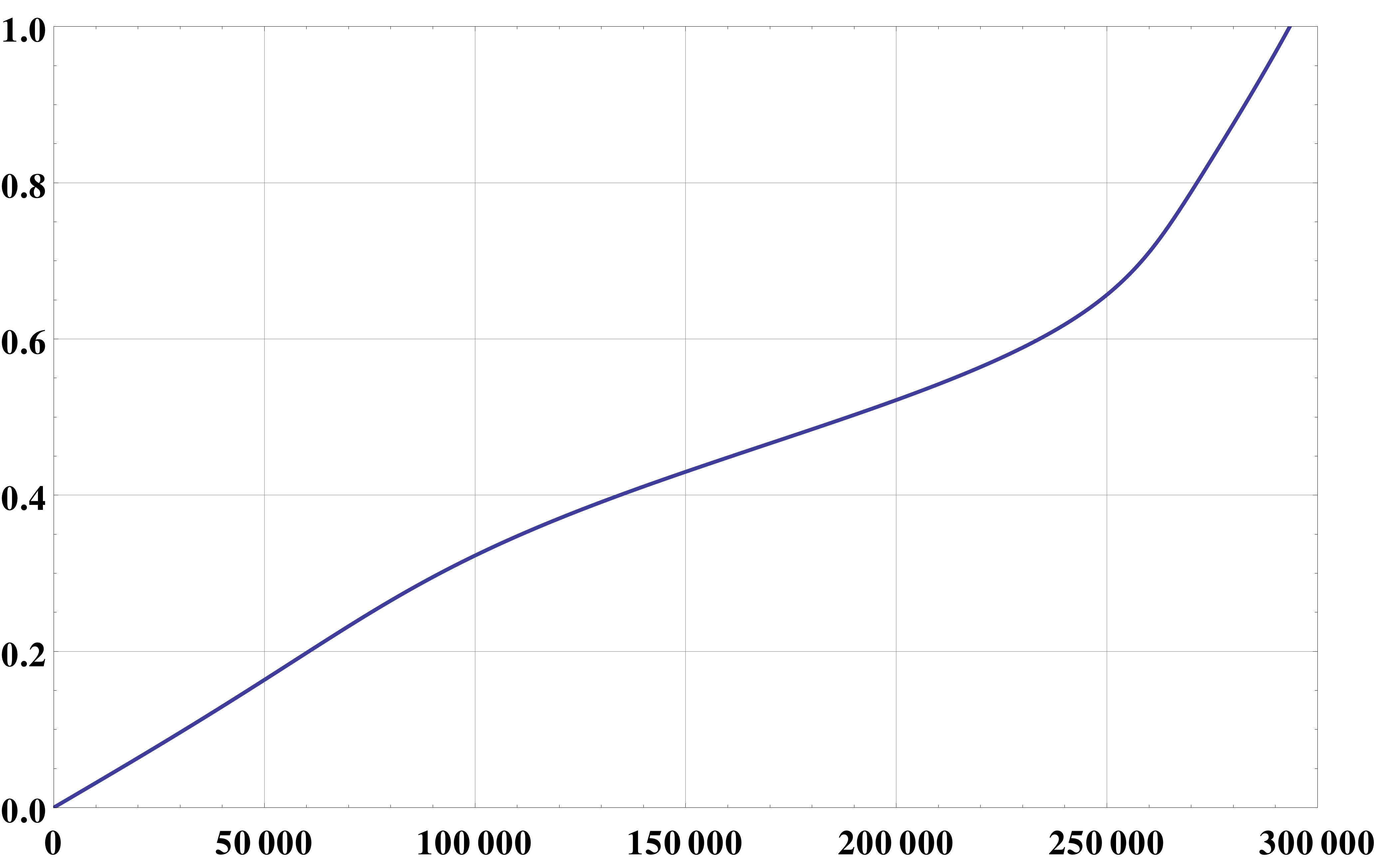}
\end{center}
\caption{\small The inverse of the nonoscillatory phase function depicted
in Figure~\ref{figure:artificial}.}
\label{figure:art2}
\end{figure}

\label{section:experiments:artificial}
\end{subsection}

%

%
%
\begin{subsection}{Gauss-Legendre quadrature formulas}

The Legendre polynomial $P_n$  of order $n$ is a solution of 
the ordinary differential equation
\begin{equation}
(1-t^2) \varphi''(t) - 2t \varphi'(t) + n(n+1) \varphi(t) = 0
\ \ \mbox{for all}\ \  -1 \leq t \leq 1.
\label{experiments:legendre:ode}
\end{equation}
Its roots 
\begin{equation}
-1 < t_1 < t_2 < \ldots < t_n < 1
\label{experiments:legendre:nodes}
\end{equation}
are the nodes of the $n$-point Gauss-Legendre 
quadrature rule, and the corresponding weights $w_1,\ldots,w_n$ are given by
\begin{equation}
w_j = \frac{2}{\left(P_n'(t_j)\right)^2 \left(1-t_j^2\right)}.
\label{experiments:legendre:weights}
\end{equation}
Formula~(\ref{experiments:legendre:weights}) can be found in many references;
it is a special case of
of Formula~15.3.1  in \cite{Szego}, for instance.
Since the Legendre polynomial
$P_n$ satisfies the symmetry relation
\begin{equation}
P_n(-t) = (-1)^n P_n(t),
\end{equation}
it suffices to compute its roots on the interval $[0,1]$.  
Rather than computing the roots of the function  $P_n$ on $[0,1]$,
we calculated the roots of the function
\begin{equation}
z_n(\theta) = P_n(\cos(\theta)) \sqrt{\sin(\theta)},
\label{experiments:legendre:transform}
\end{equation}
which is a solution of  the second order differential equation
\begin{equation}
z''(\theta) +  
\left(
\frac{1}{2} + n + n^2 + \frac{1}{4}\cot(\theta)^2
\right)
z(\theta) = 0
\label{experiments:legendre:ode2}
\end{equation}
on the interval $(0,\pi/2]$.  That (\ref{experiments:legendre:transform})
satisfies (\ref{experiments:legendre:ode2}) can be verified directly by plugging 
(\ref{experiments:legendre:transform}) into (\ref{experiments:legendre:ode2}) and
making use of the fact that $P_n$ satisfies (\ref{experiments:legendre:ode}).
     The introduction of the new dependent variable 
$\theta = \arccos(t)$ in (\ref{experiments:legendre:ode}) is suggested in  \cite{Swarztrauber}
as a way to mitigate the numerical  
problems caused by the clustering of Gauss-Legendre quadrature  
nodes near the points $\pm 1$, and the presence of the factor
$\sqrt{\sin(\theta)}$ ensures that the transformed equation is of the form
(\ref{introduction:second_order}).
See  \cite{Swarztrauber} and \cite{Bogaert-Michiels-Fostier}
for discussions of the numerical issues which arise from the clustering
of Gauss-Legendre nodes near $\pm 1$.  

\begin{table}[b!]
\small
\begin{center}
\begin{tabular}{c@{\hspace{3em}}ccccc}
\toprule
\addlinespace[.5em]
Order
&  Phase function    & Quadrature        & Total running    & Running time of   & Maximum\\
(n)&  time              & evaluation time   & time             & the algorithm     &relative difference   \\
&                    &                   &                  & of \cite{Bogaert} & in weights  \\
\midrule
\addlinespace[.75em]
$10^{3}$   & 8.70\e{-02}     & 3.74\e{-04} & 8.74\e{-02}   & 1.49\e{-04}  & 2.31\e{-14}   \\
\addlinespace[.25em]
$10^{4}$   & 9.28\e{-02}     & 2.52\e{-03} & 9.54\e{-02}   & 2.20\e{-03}  & 3.34\e{-14}   \\
\addlinespace[.25em]
$10^{5}$   & 4.45\e{-02}     & 7.37\e{-03} & 5.18\e{-02}   & 8.43\e{-03}  & 5.88\e{-14}   \\
\addlinespace[.25em]
$10^{6}$   & 4.66\e{-02}     & 7.67\e{-02} & 1.23\e{-01}   & 8.64\e{-02}  & 1.31\e{-14}   \\
\addlinespace[.25em]
$10^{7}$   & 4.62\e{-02}     & 7.21\e{-01} & 7.67\e{-01}   & 7.93\e{-01}  & 1.21\e{-14}   \\
\addlinespace[.25em]
$10^{8}$   & 4.31\e{-02}     & 7.22\e{+00} & 7.32\e{+00}   & 8.13\e{+00}  & 1.26\e{-14}   \\
\addlinespace[.25em]
$10^{9}$   & 4.79\e{-02}     & 7.02\e{+01} & 7.02\e{+01}   & 7.91\e{+01}  & 1.32\e{-14}   \\
\addlinespace[.25em]
$10^{10}$   & 4.87\e{-02}     & 7.23\e{+02} & 7.23\e{+02}   & 8.20\e{+02}  & 1.41\e{-14}   \\
\addlinespace[.25em]
\bottomrule
\end{tabular}
\caption{\small
 A comparison of the time taken to construct 
Gauss-Legendre quadrature rules of various orders
using the approach of this paper and 
 the specialized approach of \cite{Bogaert}.  There computations were performed
on a single processor core. 
All times are in seconds.}
\label{experiments:legendre:table1}
\end{center}
\end{table}

The coefficient in (\ref{experiments:legendre:ode2}) is singular
at the origin; consequently,  we represented $\alpha$ and $r$ 
using a  ``graded mesh'' of points which cluster near $0$.
More specifically, the functions $\alpha$ and $r$ were
represented via their values at the  $5$-point Chebyshev grids on the $3473$ intervals
\begin{equation}
[\gamma_1,\gamma_2],
[\gamma_2,\gamma_3], \ldots,
[\gamma_{3473},\gamma_{3474}],
\label{experiments:legendre:intervals0}
\end{equation}
where
\begin{equation}
\gamma_i = \frac{\pi}{2} \cdot (1.01)^{-3474+i}.
\end{equation}
Low order expansions were used in order to ensure that the functions
$\alpha$  and $\alpha'$ could be evaluated quickly.  
  When higher
order expansions are used, the resulting phase functions can be 
stored much more efficiently (see, for instance, the experiments
of Section~\ref{section:experiments:jacobi}).

The coefficients $c_1$ and $c_2$ such that  
\begin{equation}
z_n(\theta) = 
c_1 \frac{\cos(\alpha(\theta))}{\sqrt{\alpha'(\theta)}} +
c_2 \frac{\sin(\alpha(\theta))}{\sqrt{\alpha'(\theta)}}
\end{equation}
were obtained via Formulas~(\ref{phase:c1}) and (\ref{phase:c2});
the values of  $z_n$ and its derivative $z_n'$
at the point  $\gamma_1 = \pi/2 (1.01)^{-3473} \approx 1.54166044582463\e{-15}$,  
which are needed in (\ref{phase:c1}) and (\ref{phase:c2}), were approximated using the 
 expansion
\begin{equation}
\begin{aligned}
P_n(\cos(\theta)) \sqrt{\sin(\theta)} \approx 
\sqrt{\theta }&+\theta ^{5/2} \left(-\frac{n^2}{4}-\frac{n}{4}-\frac{1}{12}\right)+\theta ^{9/2}
   \left(\frac{n^4}{64}+\frac{n^3}{32}+\frac{5 n^2}{192}+\frac{n}{96}+\frac{1}{1440}\right)\\
&+\mathcal{O}\left(\theta^{13/2} n^6 \right).
\end{aligned}
\label{experiments:legendre:taylor1}
\end{equation}
%
%
 The real numbers $d_1$ and $d_2$ such that
\begin{equation}
z_n(\theta)  = d_1 \frac{\sin(\alpha(\theta)+d_2)}{\sqrt{\alpha'(\theta)}}
\end{equation}
were computed from $c_1$ and $c_2$ using
(\ref{phase:d1}) and (\ref{phase:d2}).  
The roots (\ref{experiments:legendre:nodes}) of $P_n$ are related to the roots 
$\theta_1 < \theta_2 < \ldots < \theta_{\lceil n/2\rceil }$
 of $z_n$ via the formula
\begin{equation}
t_j = 
\begin{cases}
-\cos(\theta_j)     &\ \ \mbox{if}\ \ 1 \leq j \leq \left\lceil \frac{n}{2}  \right\rceil\\
\cos(\theta_{n-j+1}) &\ \ \mbox{if}\ \ \left\lceil \frac{n}{2} \right\rceil < j \leq n.\\
\end{cases}
\label{experiments:legendre:nodes2}
\end{equation}
Moreover, if $t_j = \pm \cos(\theta)$, then
\begin{equation}
\frac{z_n'(\theta)}{\sqrt{\sin(\theta)}} = 
- P_n'(\cos(\theta)) \sin(\theta)^2
= 
P_n'(t_j) \sqrt{1-t_j^2}.
\label{experiments:legendre:zthetaj}
\end{equation}
By combining (\ref{phase:derivative_formula}),
(\ref{experiments:legendre:weights}), (\ref{experiments:legendre:nodes2})
and (\ref{experiments:legendre:zthetaj}), 
we obtain the formula
\begin{equation}
w_j =   
\begin{dcases}
\frac{2}{d_1^2} \frac{\sin(\theta_j)}{\alpha'(\theta_j)}    
&\ \ \mbox{if}\ \ 1 \leq j \leq \left\lceil \frac{n}{2}  \right\rceil\\
\frac{2}{d_1^2} \frac{\sin(\theta_{n-j+1})}{\alpha'(\theta_{n-j+1})}    
&\ \ \mbox{if}\ \ \left\lceil \frac{n}{2} \right\rceil < j \leq n,\\
\end{dcases}
\label{experiments:legendre:weights2}
\end{equation}
which expresses the  $j^{th}$ Gauss-Legendre weight in terms of the
derivative of the phase function $\alpha$.  
%

For several values of $n$, we compared the time 
taken to compute the nodes and
weights of the $n$-point Gauss-Legendre quadrature
via by the algorithm of  this paper with the time required to do
so via the algorithm of \cite{Bogaert}. 
We used the C++ implementation
\cite{BogaertCode}  made available by the author of \cite{Bogaert}.
These calculations were performed on a single processor core.
Table~\ref{experiments:legendre:table1} reports the results as well as
 the largest relative difference in the weights for each value of $n$.
We observe that (surprisingly, given its great generality)
 the algorithm of this paper is competitive
with that of \cite{Bogaert} when $n$ is sufficiently large.

\label{section:experiments:legendre}
\end{subsection}

%
%
\begin{subsection}{Gauss-Jacobi quadrature formulas}

 The Jacobi polynomial
$P_n^{(\gamma,\zeta)}$ is a solution of the  second order linear ordinary differential equation
\begin{equation}
(1-t^2) \varphi''(t) + \left(\zeta - \gamma- (\gamma + \zeta +2) x \right) \varphi'(t) 
+ n(n+\gamma+ \zeta + 1) \varphi(t) = 0.
\label{experiments:jacobi:ode1}
\end{equation}
Its roots
\begin{equation}
-1 < t_1 < t_2 <\ldots  < t_n < 1
\label{experiments:jacobi:nodes}
\end{equation}
are the nodes of the Gauss-Jacobi quadrature rule 
\begin{equation}
\int_{-1}^1 f(t) (1-t)^\gamma (1+t)^\zeta\ dt \approx
\sum_{j=1}^n f(t_j) w_j,
\label{experiments:jacobi:rule}
\end{equation}
which is exact for polynomials of degree $2n-1$.  The weights 
$w_1,\ldots,w_n$ are given by the formula
\begin{equation}
w_j = 
\frac{\Gamma(n+\gamma+1)\Gamma(n+\zeta+1)}{\Gamma(n+1)\Gamma(n+\gamma+\zeta+1)}
\frac{2^{\gamma+\zeta+1}}{(1-t_j^2) \left(P_n^{(\gamma,\zeta)'}(t_j)\right)^2}.
\label{experiments:jacobi:weights1}
\end{equation}
See, for instance, Section~15.3 in \cite{Szego} for more information
regarding Gauss-Jacobi quadratures, including a derivation
of Formula~(\ref{experiments:jacobi:weights1}).

\begin{table}[b!!!!!!!]
\begin{center}
\small
\begin{tabular}{ccc@{\hspace{1em}}cc@{\hspace{3em}}ccc@{\hspace{1em}}cc}
\toprule
\addlinespace[.5em]
\addlinespace[.50em]
\begin{multirow}{2}{*}{$\gamma$}\end{multirow}                &
\begin{multirow}{2}{*}{$\zeta$}\end{multirow}                &
\begin{multirow}{2}{*}{n}\end{multirow}                  &
Number of                   &   Expansion               &  
\begin{multirow}{2}{*}{$\gamma$}\end{multirow}                &
\begin{multirow}{2}{*}{$\zeta$}\end{multirow}                &
\begin{multirow}{2}{*}{n}\end{multirow}                  &
Number of                   &   Expansion                 \\
 & &   &  subintervals           &  size              & 
 & &   &  subintervals          &  size    \\
\midrule
\addlinespace[.75em]
$-0.30$ & $0.25$ & $10^{3}$    & $97$ & $2814$ & 
$0.25$ & $-0.30$ & $10^{3}$    & $97$ & $2814$ \\
\addlinespace[.25em]
 &  &$10^{4}$    & $92$ & $2669$ & 
 &  &$10^{4}$    & $92$ & $2669$ \\
\addlinespace[.25em]
 &  &$10^{5}$    & $89$ & $2582$ & 
 &  &$10^{5}$    & $89$ & $2582$ \\
\addlinespace[.25em]
 &  &$10^{6}$    & $86$ & $2495$ & 
 &  &$10^{6}$    & $86$ & $2495$ \\
\addlinespace[.25em]
 &  &$10^{7}$    & $82$ & $2379$ & 
 &  &$10^{7}$    & $82$ & $2379$ \\
\addlinespace[.25em]
 &  &$10^{8}$    & $79$ & $2292$ & 
 &  &$10^{8}$    & $79$ & $2292$ \\
\addlinespace[.25em]
 &  &$10^{9}$    & $76$ & $2205$ & 
 &  &$10^{9}$    & $76$ & $2205$ \\
\addlinespace[.25em]
 &  &$10^{10}$    & $73$ & $2118$ & 
 &  &$10^{10}$    & $73$ & $2118$ \\
\addlinespace[.25em]
 &  &$10^{11}$    & $70$ & $2031$ & 
 &  &$10^{11}$    & $70$ & $2031$ \\
\addlinespace[.25em]
 &  &$10^{12}$    & $67$ & $1944$ & 
 &  &$10^{12}$    & $67$ & $1944$ \\
\addlinespace[1em]
$\pi/2$ & $\sqrt{2}$ & $10^{3}$   & $98$ & $2843$ & 
$\sqrt{2}$ & $\pi/2$ & $10^{3}$   & $98$ & $2843$ \\
\addlinespace[.25em]
 &  &$10^{4}$    & $94$ & $2727$ & 
 &  &$10^{4}$    & $94$ & $2727$ \\
\addlinespace[.25em]
 &  &$10^{5}$    & $91$ & $2640$ & 
 &  &$10^{5}$    & $91$ & $2640$ \\
\addlinespace[.25em]
 &  &$10^{6}$    & $87$ & $2524$ & 
 &  &$10^{6}$    & $88$ & $2524$ \\
\addlinespace[.25em]
 &  &$10^{7}$    & $84$ & $2437$ & 
 &  &$10^{7}$    & $84$ & $2437$ \\
\addlinespace[.25em]
 &  &$10^{8}$    & $81$ & $2350$ & 
 &  &$10^{8}$    & $81$ & $2350$ \\
\addlinespace[.25em]
 &  &$10^{9}$    & $77$ & $2234$ & 
 &  &$10^{9}$    & $78$ & $2234$ \\
\addlinespace[.25em]
 &  &$10^{10}$    & $90$ & $2611$ & 
 &  &$10^{10}$    & $90$ & $2611$ \\
\addlinespace[.25em]
 &  &$10^{11}$    & $72$ & $2089$ & 
 &  &$10^{11}$    & $72$ & $2089$ \\
\addlinespace[.25em]
 &  &$10^{12}$    & $103$ & $2988$ & 
 &  &$10^{12}$    & $102$ & $2959$ \\
\addlinespace[.25em]
\bottomrule
\end{tabular}
\caption{\small
The size of the piecewise Chebyshev expansions used to represent
the nonoscillatory phase function representing the solution
of $z_n^{(\gamma,\zeta)}$ of Equation~(\ref{experiments:jacobi:ode2}).  The zeros
of $z_n^{(\gamma,\zeta)}$ on the interval $[0,\pi/2]$ are related to those
of the Jacobi polynomial $P_n^{(\gamma,\zeta)}$ on the interval $[0,1]$
through Formula~(\ref{experiments:jacobi:nodes2}). }
\label{experiments:jacobi:table1}
\end{center}
\end{table}

As in the special case of Gauss-Legendre quadrature rules,
we introduce a change of variables in order to avoid the problems
associated with the clustering of the nodes
(\ref{experiments:jacobi:nodes}) near $\pm 1$; see
\cite{Hale-Townsend} for further discussion of this
phenomenon.
If  $Q_n^{(\gamma,\zeta)}$ and $r^{(\gamma,\zeta)}$ are the functions defined 
via 
\begin{equation}
r^{(\gamma,\zeta)}(\theta) = 
\left(\cot\left(\frac{\theta}{2}\right)\right)^{\frac{\zeta-\gamma}{2}}
\left(\sin\left(\frac{\theta}{2}\right)\right)^{\frac{1+\gamma+\zeta}{2}}
\left(\sin(\theta)\right)^{\frac{\gamma+\zeta+1}{2}}
\label{experiments:jacobi:definition_of_r}
\end{equation}
and
\begin{equation}
\begin{aligned}
Q^{(\gamma,\zeta)}_n(\theta) = 
n (\gamma+\zeta+ &n+1)
-\frac{1}{4} \csc ^2(t) 
(\gamma -\zeta +(\gamma +\zeta +1) \cos (t))^2
\\
&+\frac{1}{2} \csc ^2(t) (\gamma+\zeta
   +(\gamma-\zeta) \cos (t)+1),
\end{aligned}
\label{experiments:jacobi:definition_of_q}
\end{equation}
then  the function $z_n^{(\gamma,\zeta)}$ defined by
\begin{equation}
z_{n}^{(\gamma,\zeta)}(\theta) = P_n^{(\gamma,\zeta)}(\cos(\theta))   r^{(\gamma,\zeta)}(\theta)
\end{equation}
satisfies the second order equation
\begin{equation}
y''(\theta)  + Q^{(\gamma,\zeta)}_n(\theta) y(\theta) = 0
\label{experiments:jacobi:ode2}
\end{equation}
on the interval $(0,\pi/2]$.  

\begin{table}[b!]
\begin{center}
\small
\begin{tabular}{c@{\hspace{2em}}ccc@{\hspace{2em}}ccc}
\toprule
\addlinespace[.5em]
& \multicolumn{3}{c}{\large $\gamma = -0.3$, \hskip 1em$\zeta = 0.25$}
& \multicolumn{3}{c}{\large $\gamma = \pi/2$, \hskip 1em$\zeta = \sqrt{2}$}\\
\addlinespace[.50em]
Order 
 &  Phase function           &  Quadrature              & Maximum 
 &  Phase function           &  Quadrature              & Maximum   \\
(n) &  time                     &  time                    &relative error  
 &  time                     &  time                    &relative error  \\ 
& & & in weights&  & & in weights\\
\midrule
\addlinespace[.75em]
$10^{3}$    & 5.14\e{-02}     & 1.02\e{-02} & 8.49\e{-14} 
           & 4.98\e{-02}     & 3.55\e{-02} & 3.59\e{-14} \\
\addlinespace[.25em]
$10^{4}$    & 4.75\e{-02}     & 3.31\e{-02} & 8.19\e{-14} 
           & 6.59\e{-02}     & 5.07\e{-02} & 4.01\e{-14} \\
\addlinespace[.25em]
$10^{5}$    & 4.62\e{-02}     & 3.85\e{-02} & 2.07\e{-14} 
           & 4.70\e{-02}     & 2.51\e{-02} & 1.43\e{-14} \\
\addlinespace[.25em]
$10^{6}$    & 4.41\e{-02}     & 3.92\e{-02} & 3.64\e{-14} 
           & 4.37\e{-02}     & 3.91\e{-02} & 2.24\e{-14} \\
\addlinespace[.25em]
$10^{7}$    & 4.13\e{-02}     & 2.17\e{-01} & 5.21\e{-14} 
           & 4.10\e{-02}     & 2.59\e{-01} & 3.68\e{-14} \\
\addlinespace[.25em]
$10^{8}$    & 4.22\e{-02}     & 2.07\e{+00} & 5.87\e{-15} 
           & 4.29\e{-02}     & 2.01\e{+00} & 1.76\e{-14} \\
\addlinespace[.25em]
$10^{9} $   & 3.95\e{-02}     & 1.82\e{+01} & 3.99\e{-15} 
           & 3.95\e{-02}     & 1.99\e{+01} & 3.52\e{-14} \\
\addlinespace[.25em]
$10^{10}$   & 4.24\e{-02}     & 1.85\e{+02} & 4.28\e{-15} 
           & 5.07\e{-02}     & 1.85\e{+02} & 1.10\e{-15} \\
\addlinespace[.25em]
$10^{11}$   & 3.77\e{-02}     & 1.76\e{+03} & 6.57\e{-15} 
           & 3.81\e{-02}     & 1.84\e{+03} & 1.29\e{-14} \\
\addlinespace[.25em]
$10^{12}$   & 3.65\e{-02}     & 1.78\e{+04} & 4.54\e{-15} 
           & 6.39\e{-02}     & 1.87\e{+04} & 9.99\e{-15} \\
\addlinespace[.25em]
\bottomrule
\end{tabular}
\caption{\small
The time  taken to compute Gauss-Jacobi quadrature rules 
of various orders via
the algorithm of this paper, and the accuracy of the resulting rules.
All times are in seconds.  
A maximum of $28$ simultaneous threads of executions were allowed during
these calculations.}
\label{experiments:jacobi:table2}
\end{center}
\end{table}

The behavior of the  coefficient in (\ref{experiments:jacobi:ode2}) 
depends strongly on  $\gamma$ and $\zeta$ --- indeed, it is singular
 in some cases and smooth on the whole interval $[0,\pi/2]$ in others.
Accordingly,  for various values of $\gamma$ and $\zeta$,
we used the adaptive version of the algorithm of Section~\ref{section:algorithm}
in order to construct a nonoscillatory phase function $\alpha^{(\gamma,\zeta)}$
 representing $z_n^{(\gamma,\zeta)}$  on the interval $(10^{-15},1)$.
Table~\ref{experiments:jacobi:table1} reports the size 
of the piecewise Chebyshev expansions used to represent
the nonoscillatory phase function in each case.
More specifically, the associated nonoscillatory phase
functions were represented via their values
at the nodes of the $30$ point Chebyshev grids on 
a collection of subintervals and Table~\ref{experiments:jacobi:table1}
lists the number of subintervals in each case and the total
number of values used to represent each of the nonoscillatory
phase functions.  We refer to this last quantity, which is equal
to $29m+1$, where $m$ is the number of subintervals into which
$[0,\pi/2]$ is divided, as the ``expansion size'' for want of a better
term.

The values of $z_n^{(\gamma,\zeta)}$ and its derivative at the point
$1.0\e{-15}$, which are needed to calculate the constants $d_1^{(\gamma,\zeta)}$ and
$d_2^{(\gamma,\zeta)}$ such that
\begin{equation}
z_n^{(\gamma,\zeta)}(\theta) = d_1^{(\gamma,\zeta)} \frac{\sin \left(\alpha^{(\gamma,\zeta)}(t) + 
d_2^{(\gamma,\zeta)}\right)}{\sqrt{\alpha^{(\gamma,\zeta)}{}'(t)}},
\end{equation}
were computed using a $7$-term Taylor expansion for $z_n^{(\gamma,\zeta)}$ 
 around the point $0$.
This expression is too cumbersome to reproduce here, but it can be
derived easily starting from the well-known representation
of $P^{(\gamma,\zeta)}_n$ in terms of Gauss' hypergeometric function
(see, for instance, Formula~(4.21.2) in \cite{Szego}).

Since
\begin{equation}
P_n^{(\gamma,\zeta)}(t) = (-1)^n P_n^{(\zeta,\gamma)}(t),
\end{equation}
 the roots of $P_n^{(\gamma,\zeta)}$  in the interval $[-1,1]$ can be obtained by computing
the roots of $z^{(\gamma,\zeta)}_n$ in the interval $(0,\pi/2]$ and
those of $z^{(\zeta,\gamma)}_n$ in the interval $(0,\pi/2)$.
More specifically, if we denote by 
\begin{equation}
\theta_1 < \theta_2 < \ldots < \theta_{\lceil n/2 \rceil}
\label{experiments:jacobi:roots1}
\end{equation}
the roots of the function $z_n^{(\zeta,\gamma)}$ in the interval $(0,\pi/2]$
and by 
\begin{equation}
\theta_{\lceil n/2 \rceil+1} < \theta_{\lceil n/2 \rceil+2}  < \ldots < \theta_{n}
\label{experiments:jacobi:roots2}
\end{equation}
the roots of $z_n^{(\gamma,\zeta)}$ in the interval  $(0,\pi/2)$, then
the nodes of the $n$-point Gauss-Jacobi quadrature rule are given by the
formula
\begin{equation}
t_j = 
\begin{cases}
-\cos(\theta_j)     &\ \ \mbox{if}\ \ 1 \leq j \leq \left\lceil \frac{n}{2}  \right\rceil\\
\cos(\theta_{n-j+1}) &\ \ \mbox{if}\ \ \left\lceil \frac{n}{2} \right\rceil < j \leq n.\\
\end{cases}
\label{experiments:jacobi:nodes2}
\end{equation}
%
As in the case of Gauss-Legendre quadrature rules,
the weights of Gauss-Jacobi quadrature rules
can be  expressed  in terms
of the derivatives of the phase functions which represent the functions
$z_n^{(\gamma,\zeta)}$ and $z_n^{(\zeta,\gamma)}$;
more specifically, we combine
(\ref{phase:derivative_formula}), (\ref{experiments:jacobi:weights1})
and (\ref{experiments:jacobi:nodes2}) to obtain
\begin{equation}
w_j = \begin{dcases}
\frac{\Gamma(n+\gamma+1)\Gamma(n+\zeta+1)}{\Gamma(n+1)\Gamma(n+\gamma+\zeta+1)}
\frac{\left(r^{(\zeta,\gamma)}(\theta_j)\right)^2}{\left(d_1^{(\zeta,\gamma)}\right)^2 \alpha^{(\zeta,\gamma)}{}'(\theta_j)}
    &\ \ \mbox{if}\ \ 1 \leq j \leq \left\lceil \frac{n}{2}  \right\rceil\\
\frac{\Gamma(n+\gamma+1)\Gamma(n+\zeta+1)}{\Gamma(n+1)\Gamma(n+\gamma+\zeta+1)}
\frac{\left(r^{(\gamma,\zeta)}(\theta_{n-j+1})\right)^2}{\left(d_1^{(\gamma,\zeta)}\right)^2 \alpha^{(\gamma,\zeta)}{}'
(\theta_{n-j+1})}
    &\ \ \mbox{if}\ \ \left\lceil \frac{n}{2}  \right\rceil< j \leq n.
\end{dcases}
\label{experiments:jacobi:weights2}
\end{equation}
 We follow \cite{Hale-Townsend} in using the asymptotic formula
\begin{equation}
\frac{\Gamma(\gamma+n) \Gamma(\zeta+n)}{\Gamma(\chi+n)\Gamma(\gamma+\zeta-\chi+n)}
\approx
1 + \sum_{m=1}^M \frac{(-\gamma)_m(-\zeta)_m}{\Gamma(m+1) (-\gamma-\zeta-n)_m},
\label{experiments:jacobi:gamma_approximation}
\end{equation}
which is a special case of (3.1) in \cite{Buhring}, in order to 
evaluate the ratio of gamma functions
appearing in (\ref{experiments:jacobi:weights2}).
The symbol $(x)_m$ appearing in (\ref{experiments:jacobi:gamma_approximation})
is the Pochhammer symbol, which is defined via
\begin{equation}
(x)_m = x (x+1) \ldots (x+m).
\end{equation}

We used the algorithm of Section~\ref{section:algorithm} to construct Gauss-Jacobi rules
of various orders $n$ and for various values of $\gamma$ and $\zeta$.
We tested the accuracy of these rules by comparing
the first $\min\{10^7,n\}$  weights
to those generated by running the Glaser-Liu-Rokhlin \cite{Glaser-Rokhlin} algorithm
using IEEE quadruple precision arithmetic.
Table~\ref{experiments:jacobi:table2} reports the  results.
For each combination of $\gamma$, $\zeta$ and $n$ considered,
it lists the time taken to compute the nonoscillatory phase
functions representing $z_n^{(\gamma,\zeta)}$ and $z_n^{(\gamma,\zeta)}$ and their inverses, the total time
required to calculate the nodes and weights of 
the corresponding Gauss-Jacobi quadrature rule,
and the maximum relative error in the weights of that Gauss-Jacobi rule.
A maximum of $28$ simultaneous threads of execution were allowed
during these calculations.  

In Table~\ref{experiments:jacobi:table3}, we compare the 
time required to compute Gauss-Jacobi quadratures rules of various
orders using the algorithm of this paper with the time required
to do so 
using the algorithm of \cite{Hale-Townsend}.   The parameters
were taken to be $\gamma = 0.2$ and $\zeta = 0.5$.  We used the adaptive
version of our algorithm with $k = 30$ and 
the Julia implementation \cite{TownsendCode} provided by the authors
of \cite{Hale-Townsend}.
These calculations were performed on a single processor core.  

\begin{table}[h!]
\small
\begin{center}
\begin{tabular}{c@{\hspace{3em}}cc}
\toprule
\addlinespace[.5em]
Order & Algorithm of     & Algorithm of         \\
      & this paper       & \cite{Hale-Townsend} \\
\midrule
$10^3$ & 5.58 \e{-02} & 3.22 \e{-02} \\
\addlinespace[.25em]
$10^4$ & 5.50 \e{-02} & 1.65 \e{-01} \\
\addlinespace[.25em]
$10^5$ & 8.63 \e{-02} & 1.29 \e{+00} \\
\addlinespace[.25em]
$10^6$ & 4.72 \e{-01} & 1.14 \e{+01} \\
\addlinespace[.25em]
$10^7$ & 3.13 \e{+00} & 2.03 \e{+02} \\
\addlinespace[.25em]
$10^8$ & 2.92 \e{+01} & 2.03 \e{+03} \\
\addlinespace[.25em]
$10^9$ & 2.89 \e{+02} & --- \\
\addlinespace[.25em]
$10^{10}$ & 2.82 \e{+03} & --- \\
\bottomrule
\end{tabular}
\caption{\small
 A comparison of the time taken to construct 
Gauss-Jacobi quadrature rules 
of various orders 
using the approach of this paper and 
 the specialized approach of \cite{Hale-Townsend}.  
Here, the parameters were taken to be
$\gamma = 0.2$ and $\zeta = 0.5$.
There computations were performed
on a single processor core and all times are in seconds.
Entries marked with a ``---'' indicate experiments which were prohibitively expensive
to perform.}
\label{experiments:jacobi:table3}
\end{center}
\end{table}

\label{section:experiments:jacobi}
\end{subsection}


%
%
\begin{subsection}{Gauss-Laguerre quadrature formulas}

 The Laguerre polynomial $L_n^{(\gamma)}$ 
 is a solution of the ordinary differential equation
\begin{equation}
t \psi''(t) + (1+\gamma-t) \psi'(t) + n \psi(t) = 0 
\ \ \mbox{for all} \ \ 0 \leq t < \infty.
\label{experiments:laguerre:ode1}
\end{equation}
Its zeros $t_1 < t_2 < \ldots < t_n$
are the nodes of the Gauss-Laguerre quadrature rule
\begin{equation}
\int_{0}^\infty t^\gamma\exp(-t) f(t)\ dt \approx \sum_{k=1}^n f(t_k) w_k,
\label{experiments:laguerre:rule}
\end{equation}
and the weights $w_1,\ldots,w_n$ are given by the formula
\begin{equation}
w_j = \frac{\Gamma(n+\gamma+1)}{\Gamma(n+1)}\frac{1}{t_j \left(L_n^{(\gamma)}{}'(t_j)\right)^2}.
\label{experiments:laguerre:weights}
\end{equation}
Formula~(\ref{experiments:laguerre:weights}) can be 
found in many sources; it appears as (15.3.5) in  \cite{Szego}, for instance.
The function $z_n$ defined via 
\begin{equation}
z_n(u) = 
L_n (\exp(u)) \exp\left(-\frac{\exp(u)}{2} + \frac{\gamma u}{2} \right)
\label{experiments:laguerre:zn}
\end{equation}
satisfies the second order differential equation
\begin{equation}
z''(u) +\left(\frac{\exp(u)}{2} - \frac{1}{4} \left(\gamma - \exp(u)\right)^2 + \exp(u)n\right) z(u) = 0
\label{experiments:laguerre:ode2}
\end{equation}
on the interval $(-\infty,\infty)$ and the function
\begin{equation}
y_n(v) = L_n (v^2) \exp\left(-v^2/2\right) v^{1/2 + \gamma}
\label{experiments:laguerre:yn}
\end{equation}
is a solution of 
\begin{equation}
y''(v)   + \left(2 + 2\gamma + 4n + \frac{1-4\gamma^2}{4v^2} -v^2 \right) y(v) = 0
\label{experiments:laguerre:ode3}
\end{equation}
on the interval $(0,\infty)$.

\begin{table}[b!]
\begin{center}
\small
\begin{tabular}{c@{\hspace{2em}}cc@{\hspace{2em}}cc@{\hspace{2em}}cc}
\toprule
\addlinespace[.5em]
& \multicolumn{2}{c}{\large $\gamma = -0.5$}
& \multicolumn{2}{c}{\large $\gamma = 0$ }
& \multicolumn{2}{c}{\large $\gamma = 0.5$ }\\
\addlinespace[.50em]
Order  &  Phase function &  Total      & Phase function &  Total & Phase function  &  Total                      \\
(n)    &  expansion size &  time       & expansion size &  time  & expansion size  &  time  \\
\midrule
\addlinespace[.75em]
$10^{3}$    & $8764$  & 2.37\e{-01}  
           & $8212$  & 2.16\e{-01}  
           & $8787$  & 2.29\e{-01}  \\
\addlinespace[.25em]
$10^{4}$    & $11938$ & 3.09\e{-01}  
           & $11869$ & 3.11\e{-01}  
           & $11892$ & 3.13\e{-01}  \\
\addlinespace[.25em]
$10^{5}$    & $12168$ & 3.90\e{-01}  
           & $12168$ & 3.91\e{-01}  
           & $12099$ & 3.94\e{-01}  \\
\addlinespace[.25em]
$10^{6}$    & $8741$ & 5.28\e{-01}  
           & $8695$ & 5.10\e{-01}  
           & $8419$ & 5.02\e{-01}  \\
\addlinespace[.25em]
$10^{7}$    & $4831$ & 2.71\e{+00}  
           & $4785$ & 2.71\e{+00}  
           & $2025$ & 2.66\e{+00}  \\
\addlinespace[.25em]
$10^{8}$    & $14744$ & 2.43\e{+01}  
           & $10765$ & 2.39\e{+01}  
           & $19114$ & 2.45\e{+01}  \\
\addlinespace[.25em]
$10^{9}$    & $23275$ & 2.52\e{+02}  
           & $9802$  & 3.01\e{+02}  
           & $5743$ & 3.42\e{+02}  \\
\addlinespace[.25em]
\bottomrule
\end{tabular}
\caption{\small
The time  (in seconds) taken to compute Gauss-Laguerre quadrature rules 
of various orders via
the algorithm of this paper, and the size of the expansions
of the phase functions used to represent the solution.
These calculations were performed on a single processor core.}
\label{experiments:laguerre:table1}
\end{center}
\end{table}

For each of several values of $n$ and $\gamma$, we used the algorithm of this paper
to construct two nonoscillatory phase functions, $\alpha_1$ and 
$\alpha_2$.  The function $\alpha_1$ represented solutions of 
 (\ref{experiments:laguerre:ode2})  on the interval
\begin{equation}
\left(-30, 0.0 \right),
\label{experiments:laguerre:interval}
\end{equation}
and $\alpha_2$ represented solutions of 
(\ref{experiments:laguerre:ode3}) on the interval
\begin{equation}
\left(\gamma, \sqrt{2 n + \gamma - 2 + \sqrt{1+4(n-1)(n+\gamma-1)}}\right),
\label{experiments:laguerre:interval2}
\end{equation}
where $\zeta$ is the largest root of the Laguerre polynomial
$L_n^{(\gamma)}$ in $(0,1)$.
The interval (\ref{experiments:laguerre:interval2})
was chosen in light of  the bound
\begin{equation}
t_n <  2 n + \gamma - 2 + \sqrt{1+4(n-1)(n+\gamma-1)},
\end{equation}
which can be found in  \cite{IsmailLi}.  
The phase function $\alpha_1$  representing solutions of 
(\ref{experiments:laguerre:ode2}) was constructed first;
the values of $z_n$ and its derivative at the point $-30$,
which were used in order to obtain  $c_1$ and $c_2$ such that  
\begin{equation}
z_n(u) = 
c_1 \frac{\sin(\alpha_1(u)+c_2)}{\sqrt{\alpha_1'(u)}},
\end{equation}
were calculated using a $7$-term Taylor expansion for the function
\begin{equation}
L_n (t)  \exp\left(-\frac{t}{2}\right) t^{\gamma/2}.
\end{equation}
The left endpoint of the interval (\ref{experiments:laguerre:interval})
was chosen in lieu of a bound for the smallest root of $L_n$
(of the sort appearing in \cite{IsmailLi}) in order to ensure
the accuracy of these approximations.  Once the
coefficients $c_1$ and $c_2$ were obtained, we calculated
the location of the largest root $\zeta$ of $z_n$ on the interval
$(0,1)$ as well as the value of $z_n'(\zeta)$.  
These values were used to construct the coefficients $d_1$ and $d_2$
in the  representation
\begin{equation}
y_n(v) = 
d_1 \frac{\sin(\alpha_2(v)+d_2)}{\sqrt{\alpha_2'(v)}}
\end{equation}
of $y_n$ in terms of the phase function $\alpha_2$ for
 equation (\ref{experiments:laguerre:ode3}) on the
interval (\ref{experiments:laguerre:interval2}). 
Note that, as discussed in Section~\ref{section:phase}, there is no
need to evaluate trigonometric functions of large arguments in order
to obtain the value of $z_n'$ at $\zeta$ (and hence none of the attendant
loss of precision).
We used two phase functions to represent $L_n$ because
some precision was lost  when we represented $L_n$ on an interval containing all of its zeros
using  a single phase function.  
The adaptive version of the algorithm of Section~\ref{section:algorithm}
was used to construct both phase functions and the parameter $k$ was taken to be $30$.

For each pair of chosen values of $n$ and $\gamma$, we 
computed the zeros $u_1 < u_2 < \ldots < u_k$ of $z_n$
on the interval (\ref{experiments:laguerre:interval2})
and then used the formulas
\begin{equation}
t_j = \exp(u_j)
\end{equation}
and
\begin{equation}
w_j = 
\frac{\exp(-\exp(u_j)) \exp((1+\gamma)u_j)}{c_1^2 \alpha_1'(u_j)}
\label{experiments:laguerre:weights2}
\end{equation}
in order to construct the nodes  
of (\ref{experiments:laguerre:rule}) in the interval $(0,1)$ and the corresponding weights.
The expression (\ref{experiments:laguerre:weights2})  is obtained
by combining (\ref{phase:derivative_formula}), (\ref{experiments:laguerre:weights}) and 
(\ref{experiments:laguerre:zn}).  We next computed the zeros 
$v_1 < v_2 < \ldots < v_{n-k}$ of $y_n$ in the interior 
of the interval (\ref{experiments:laguerre:interval2}).
The nodes $t_{k+1} < t_{k+2} < \ldots < t_n$ of the rule (\ref{experiments:laguerre:rule})
contained in the interval $(1,\infty)$ are related
to those of $y_n$ via the formula $t_{k+j} = v_j^2$,
and the corresponding weights are given by
\begin{equation}
w_{k+j} =  \frac{4 \exp(-v_j^2) v_j^{1+2\gamma}}{d_1^2 \alpha'(v_j)}.
\label{experiments:laguerre:weights3}
\end{equation}
Formula~(\ref{experiments:laguerre:weights3}) is obtained in the 
usual fashion --- by combining
(\ref{phase:derivative_formula}), (\ref{experiments:laguerre:weights}) and 
 (\ref{experiments:laguerre:yn}).

Table~\ref{experiments:laguerre:table1} reports the results of these experiments;
for each chosen pair of $n$ and $\gamma$,
it lists the total time required to compute the quadrature rule (including
the time required to compute the phase function and its inverse),
and the sum of the sizes of the expansions used to represent the two
nonoscillatory phase functions.
These calculations were performed on a single
processor core.


\label{section:experiments:laguerre}
\end{subsection}

%
%

\begin{subsection}{Roots of Bessel functions}

For each positive real number $\nu$, we denote by $J_\nu$ the  solution of  Bessel's equation
\begin{equation}
t^2 y''(t) + t y'(t) + (t^2-\nu^2)y(t) = 0
\ \ \ \mbox{for all} \ \ 0 \leq t < \infty
\label{experiments:bessel:eq}
\end{equation}
which is finite at the origin.      The function $J_\nu$ has an infinite number of roots
%
in the interval $(\nu,\infty)$ on which it oscillates.
The equation (\ref{experiments:bessel:eq})
is brought into the form
\begin{equation}
z''(u)+ \left(\exp(2u)-\nu^2\right)z(u) = 0
\ \ \ \mbox{for all} \ \ -\infty < u < \infty
\label{experiments:bessel:ode2}
\end{equation}
via the transformation
\begin{equation}
z(u) = y(\exp(u)).
\label{experiments:bessel:transform}
\end{equation}
For various values of $\nu$, we constructed a nonoscillatory phase
function $\alpha$ representing solutions of (\ref{experiments:bessel:ode2})
in the interval
\begin{equation}
\left[ \log(\nu), \log\left( \left(10^9 +\frac{\nu}{2}- \frac{1}{4}\right )\pi\right)\right].
\label{experiments:bessel:interval}
\end{equation}
The right endpoint in  (\ref{experiments:bessel:interval}) is an upper bound
for the location of the one billionth root of the function $J_n(\exp(u))$
(see, for instance, 10.21.19 in \cite{NISTHandbook}).
The Equation~(\ref{experiments:bessel:ode2}) has a turning point
at $u=\log(\nu)$; consequently, any phase function representing
its solutions is singular there.  We used the adaptive
version of the algorithm of Section~\ref{section:algorithm}
in order to construct phase functions in these experiments.  The parameter
$k$ was taken to be $30$.

The coefficients $d_1$ and $d_2$ such that 
\begin{equation}
J_\nu(\exp(t)) = d_1 \frac{\sin(\alpha(t) + d_2) }{\sqrt{\alpha'(t)}}
\end{equation}
were calculated using the approximations of the values of
$J_\nu$ and its derivative at the point $\nu$ 
obtained via the formulas
\begin{equation}
J_\nu(\nu) = \frac{1}{\pi} \int_0^\pi  \exp\left(-\nu F(t) \right) \ dt
\label{experiments:bessel:jnu}
\end{equation}
and
\begin{equation}
J_\nu'(\nu) = \frac{1}{\pi} \int_0^\pi  \frac{t-\sin(t)\cos(t)}{\sqrt{t^2-\sin(t)^2}}
\exp\left(-\nu F(t) \right) \ dt,
\label{experiments:bessel:jnuprime}
\end{equation}
where
\begin{equation}
F(t) = \log\left( \frac{t + \sqrt{t^2 - \sin(t)^2}}{\sin(t)}\right)
- \cot(t) \sqrt{t^2 - \sin(t)^2}.
\end{equation}
Formulas (\ref{experiments:bessel:jnu}) and (\ref{experiments:bessel:jnuprime})
appear in  Section~8.53 of \cite{Watson}, among many other  sources.
Note that the integrands in (\ref{experiments:bessel:jnu}) and
(\ref{experiments:bessel:jnuprime}) are nonoscillatory 
so that the order of the quadrature rule
 needed to calculate them does
not depend on $\nu$.

\begin{table}[t!!!!!!!]
\begin{center}
\begin{tabular}{c@{\hspace{1em}}cccc}
\toprule
\addlinespace[.5em]
$\nu$   & Phase function           &  Phase function       & Root calculation   &  
Maximum relative            \\
        & expansion size              &  time                 & time   &  error                       \\
\midrule
\addlinespace[.75em]
$\sqrt{2} \cdot 10^3$   & 2205 & 2.89\e{-02} & 8.58\e{+00} & 1.83\e{-15} \\
\addlinespace[.25em]
$\pi \cdot 10^4$        & 2305 & 3.48\e{-02} & 8.58\e{+00} & 1.81\e{-15} \\
\addlinespace[.25em]
$\pi \cdot 10^5$        & 2466 & 3.19\e{-02} & 8.63\e{+00} & 3.89\e{-14} \\
\addlinespace[.25em]
$\sqrt{3} \cdot 10^6$   & 2066 & 3.00\e{-02} & 8.48\e{+00} & 1.59\e{-15} \\
\addlinespace[.25em]
$\pi \cdot 10^7$        & 1770 & 2.50\e{-02} & 8.78\e{+00} & 1.72\e{-15} \\
\addlinespace[.25em]
$\sqrt{2} \cdot 10^8$   & 2466 & 3.22\e{-02} & 8.53\e{+00} & 1.67\e{-15} \\
\addlinespace[.25em]
$\pi \cdot 10^9$        & 3858 & 4.99\e{-02} & 8.97\e{+00} & 4.06\e{-15} \\
\addlinespace[.25em]
$\sqrt{3} \cdot 10^{10}$ & 4148 & 5.01\e{-02} & 1.06\e{+01} & 1.65\e{-15} \\
\addlinespace[.5em]
\bottomrule
\end{tabular}
\caption{\small
The time (in seconds) taken to compute the first one billion roots
of Bessel functions of various orders, the accuracy of the obtained roots,
and the size of the piecewise expansion used to represent the associated
nonoscillatory phase functions. 
A maximum of $28$ simultaneous threads of execution were allowed during these
calculations.
}
\label{experiments:bessel:table1}
\end{center}
\end{table}

For each chosen value of $\nu$, we used the nonoscillatory phase function
to compute the first one billion roots
of $J_\nu$.  The obtained values were compared against
those generated by running the Glaser-Liu-Rokhlin algorithm \cite{Glaser-Rokhlin}
using IEEE quadruple  precision arithmetic.
Table~\ref{experiments:bessel:table1} shows the
results; it reports the number of values used to represent
each nonoscillatory phase function,
the time taken to construct each phase function,
the time required to calculate the roots, and the maximum relative error in
the obtained roots.


\label{section:experiments:bessel}
\end{subsection}

\label{section:experiments}
\end{section}

\begin{section}{Conclusions}

We have described a fast and highly accurate algorithm for the computation of the roots of special
functions satisfying second order ordinary differential equations.
 Despite its great generality,
when it is used to construct classical Gaussian quadrature rules
our algorithm is competitive with specialized, state-of-the-art methods.
It is based on two  observations:  (1) the solutions of
 second order linear ordinary differential equations of the form
\begin{equation}
y''(t) + q(t) y(t) = 0,
\label{conclusion:ode}
\end{equation}
where $q$ is smooth and positive , can be represented as
\begin{equation}
y(t) = \frac{d_1\sin(d_2 + \alpha(t))}{\sqrt{\alpha'(t)}}
\label{conclusion:y}
\end{equation}
with $\alpha$ a nonoscillatory function even when the magnitude of $q$ is large, and 
(2) the roots  of a function $y$
represented in the form (\ref{conclusion:y})
and the values of its derivative $y'$ at those roots can be calculated
without evaluating  trigonometric functions of large orders and the concomitant
loss of precision.

Our algorithm was designed for reliability and accuracy at the expense of speed.  
It can be significantly 
accelerated in many cases of interest  and 
improvements in it 
will be reported  by the author at a later date.
As will algorithms for the fast evaluation
of certain special functions and the fast application of their
associated transforms which
take advantage of the fact that
explicit formulas for nonoscillatory phase functions are
sometimes available.

\label{section:conclusion}
\end{section}

\begin{section}{Acknowledgments}
The author would like to thank Vladimir Rokhlin for
several useful discussions related to this work,
and the anonymous reviewers for there many helpful comments.
The author was supported by a fellowship from the Alfred P. Sloan
Foundation, and by  National Science Foundation grant DMS-1418723.
\end{section}

\begin{section}{References}
\bibliographystyle{acm}
\bibliography{zeros}
\end{section}

\vfill\eject

\end{document}